\documentclass{ceb-l} 

\usepackage{amsmath, amsxtra, amsthm, amsfonts, amssymb}
\usepackage{mathrsfs}
\usepackage{graphicx}
\usepackage{url}
\usepackage{color}
\usepackage{bbm}
\usepackage{tikz-cd}

\usepackage[T1]{fontenc}
\usepackage{enumitem}

 \usepackage[hidelinks,backref=page]{hyperref} 
 \hypersetup{
    colorlinks,
    citecolor=magenta,
    filecolor=magenta,
    linkcolor=blue,
    urlcolor=black
}

\usepackage{cleveref}

\newcommand{\RR}{\mathbb R}
\newcommand{\CC}{\mathbb C}
\newcommand{\PP}{\mathbb P}
\newcommand{\ZZ}{\mathbb Z}
\newcommand{\M}{\mathrm{M}}
\newcommand{\kk}{\mathbbm{k}}
\newcommand{\rk}{\operatorname{rk}}
\newcommand{\U}{\mathrm{U}}
\newcommand{\T}{\mathrm{T}}
\newcommand{\be}{\mathbf{e}}

\theoremstyle{definition}
\newtheorem{thm}{Theorem}[section]
\newtheorem{cor}[thm]{Corollary}

\newtheorem{prop}[thm]{Proposition}
\newtheorem{defn}[thm]{Definition}

\newtheorem{eg}[thm]{Example}
\newtheorem{rem}[thm]{Remark}

\newtheorem{exer}[thm]{Exercise}

\title[Combinatorial essence of independence]{Essence of independence:\\ 
Hodge theory of matroids since June Huh}
\author{Christopher Eur}
\date{}

\address{Harvard University}
\email{ceur@math.harvard.edu}

\subjclass[2020]{05B35, 05E14, 14F43, 14C17}

\begin{document}

\maketitle

\begin{abstract}
Matroids are combinatorial abstractions of independence, a ubiquitous notion that pervades many branches of mathematics.  June Huh and his collaborators recently made spectacular breakthroughs by developing a Hodge theory of matroids that resolved several long-standing conjectures in matroid theory.
We survey the main results in this development and ideas behind them.
\end{abstract}


\section{Introduction}

The notion of ``independence'' resides everywhere, for example in graphs, vector configurations, field extensions, hyperplane arrangements, matchings, and discrete optimizations.
Matroid theory captures the combinatorial essence of ``independence'' shared in these structures.
For example, let us consider the following graph $G$ with edges labelled $\{1,\ldots, 5\}$ and the set of vectors $\{v_1, \ldots, v_5\}$.
\begin{figure}[h]
	\centering
\begin{tabular}{ccc}

\raisebox{10pt}{
\begin{tikzpicture}[scale=0.8]
\draw[fill]
(0,1.1) circle [radius=0.1] 
(0,-1.1) circle [radius=0.1] 
(1.73, 0) circle [radius=0.1] 
(-1.73,0) circle [radius=0.1] 
;
\draw
(0,1.1) -- (0,-1.1) node [left] [midway] {$3$}
(1.73, 0) -- (0,1.1) node [above] [midway] {$4$}
(0,-1.1) -- (1.73, 0) node [below] [midway] {$5$}
(-1.73,0) -- (0,1.1) node [midway] [above] {$1$}
(-1.73,0) -- (0,-1.1) node [midway] [below] {$2$}
;
\end{tikzpicture}
}

&\qquad
\includegraphics[height=31mm]{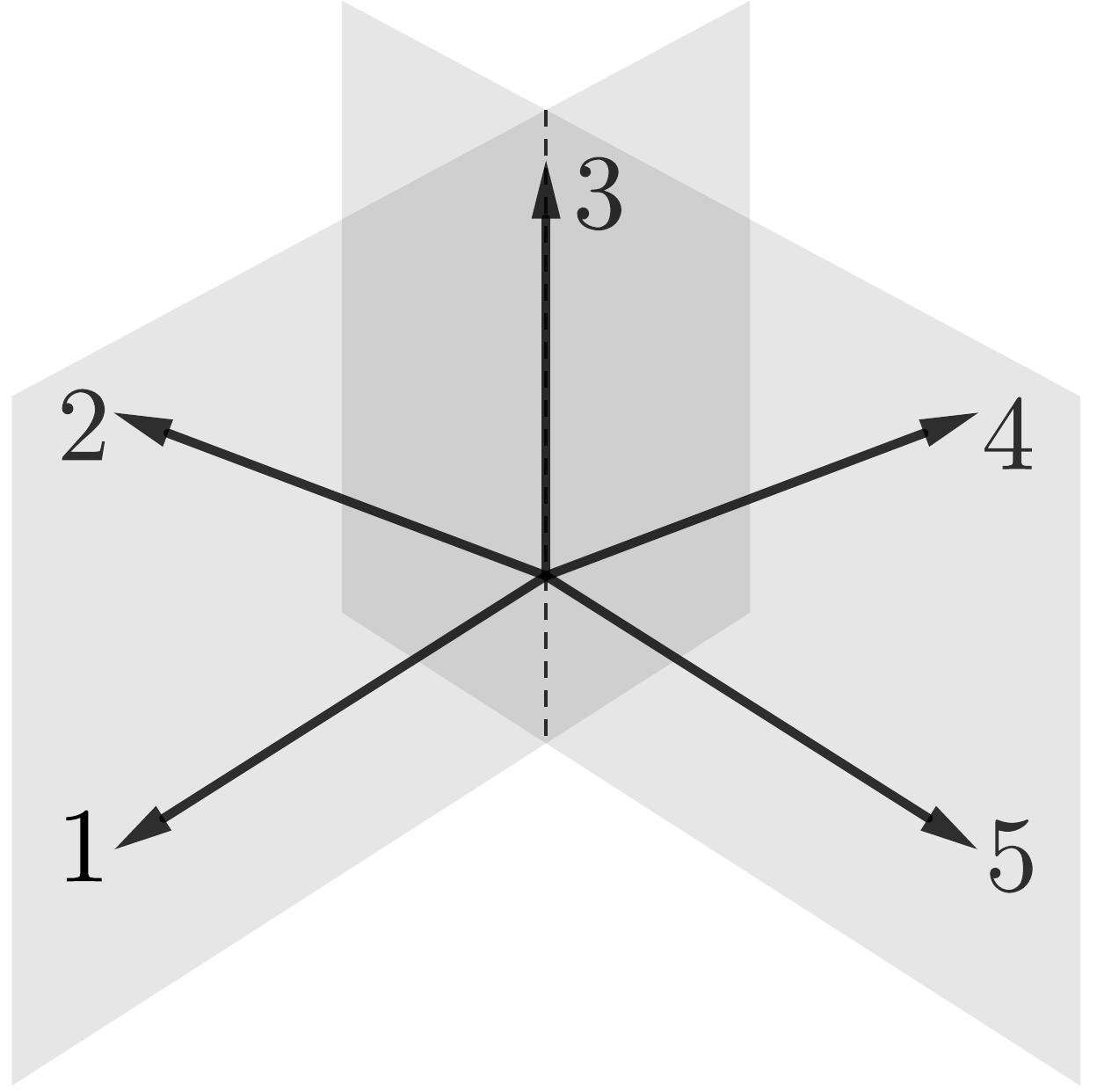}

&\quad
\raisebox{40pt}{
$
\begin{matrix}
v_1  \  \ v_2  \ \ v_3  \ \ v_4  \ \ v_5 \\
\begin{bmatrix}
1 & 1 & 0 & 0 & 0\\
0 & 1 & 1 & 1 & 0\\
0 & 0 & 0 & 1 & 1\\
\end{bmatrix}
\end{matrix}
$
}
\end{tabular}
\caption{}\label{fig:intro}
\end{figure}

\noindent We observe a common combinatorial structure:  a subset of edges in $G$ is acyclic if and only if the corresponding subset of vectors is linearly independent.
This combinatorial structure is encoded as a matroid, introduced by Whitney \cite{Whi32}.

\begin{defn}
A \textbf{matroid} $\M = (E, \mathcal I)$ consists of a finite set $E = \{1, \ldots, n\}$, called its \textbf{ground set}, and a nonempty collection $\mathcal I$ of subsets of $E$, called the \textbf{independent sets} of $\M$, such that
\begin{itemize}
\item[] if $I\in \mathcal I$ and $J\subseteq I$, then $J\in \mathcal I$, \quad and
\item[] if $I, J \in \mathcal I$ and $|I|<|J|$, then there exists an element $j\in J\setminus I$ such that $I \cup \{j\} \in \mathcal I$.
\end{itemize}
The definition implies that every maximal independent set of $\M$ has the same cardinality $r$, which we call the \textbf{rank} of $\M$.
\end{defn}

Graphs and vector spaces give prototypical examples of matroids:
\begin{itemize}
\item When $E$ is identified with the set of edges of a finite graph $G$, setting 
\[
\mathcal I = \{I \subseteq E : \text{the subset $I$ of edges in $G$ is acyclic}\}
\]
defines a matroid $\M = (E,\mathcal  I)$.  Matroids arising in this way are called \textbf{graphical matroids}.

\item When $E$ is identified with a finite set of vectors spanning a vector space $V$, setting
\[
\qquad \mathcal I = \{I \subseteq E : \text{the subset $I$ of vectors in $V$ is linearly independent}\}
\]
defines a matroid $\M = (E,\mathcal  I)$.
Matroids arising in this way are called \textbf{realizable matroids}.
\end{itemize}
We see that the graph $G$ and the set of vectors in \Cref{fig:intro} define the same matroid.

\subsection{Combinatorial sequences from a matroid}

Several long-standing conjectures in matroid theory, recently resolved by June Huh and his collaborators, concern the behavior of sequences of invariants of a matroid.  For a sequence $(a_0, a_1, \ldots, a_m)$ of nonnegative real numbers, we say that it
\begin{itemize}
\item is \emph{unimodal} if there exists $0\leq k \leq m$ such that
\[
a_0 \leq a_1 \leq \cdots \leq a_k \geq a_{k+1} \geq \cdots \geq a_m,
\]
\item is \emph{log-concave} if $a_i^2 \geq a_{i-1}a_{i+1}$ for all $1\leq i \leq m-1$,
\item has \emph{no internal zeros} if $a_ia_j \neq 0$ implies $a_k \neq 0$ for all $0\leq i \leq k \leq j \leq 0$, and
\item is \emph{top-heavy} if $a_i \leq a_{d-i}$ for all $0\leq i \leq \frac{d}{2}$ where $d$ is the largest index such that $a_d \neq 0$.
\end{itemize}
Note that a log-concave sequence is unimodal if and only if it has no internal zeroes.  For a survey of unimodality and log-concavity in combinatorics, see \cite{Sta89, Bre94}.

\smallskip
We consider the following sequences of invariants of a rank $r$ matroid $\M$.
For some of them, we describe them only for a graphical or a realizable matroid, postponing their descriptions for arbitrary matroids to \Cref{sec:matroids}.

\begin{enumerate}[label=(\alph*)]
\item \label{seqI} For $0\leq i \leq r$, let $I_i$ be the number of independent sets of $\M$ of cardinality $i$.  In other words, the sequence $(I_0, \ldots, I_r)$ is the $f$-vector of the simplicial complex whose faces are the independent sets of $\M$.

\item \label{seqI'} We may also consider the $h$-vector.  That is, for $0\leq i \leq r$, let $I'_i$ be defined by the identity $\sum_{i=0}^r I'_i q^{r-i} = \sum_{i=0}^r I_i (q-1)^{r-i}$, where $I_i$ is as in \ref{seqI}.

\item \label{seqw} Suppose $\M$ is the graphical matroid of a finite connected nontrivial graph $G$.  The \textbf{chromatic polynomial} $\chi_G(q)$ of $G$
is defined as
\[
\chi_G(q) = \text{the number of proper colorings of $G$ with at most $q$ colors},
\]
where a coloring of the vertices is \emph{proper} if no edge has both vertices the same color.  It is polynomial in $q$ of degree $r+1$, and is always divisible by $q(q-1)$.
Let $(w_0, \ldots, w_r)$ be the absolute values of the coefficients of $\frac{1}{q}\chi_G(q)$, starting from the highest degree term.

\item \label{seqw'} Continuing the assumption that $\M$ is the graphical matroid of $G$, we define $(w'_0, \ldots, w'_{r-1})$ as the absolute values of the coefficients of $\frac{1}{q(q+1)}\chi_G(q+1)$.

\item \label{seqW} Suppose $\M$ is the realizable matroid of a set of vectors $\{v_1, \ldots, v_n\}$ spanning a vector space $V$.  Let $(W_0, \ldots, W_r)$ be a sequence defined by setting for each $0\leq i \leq r$,
\[
\begin{matrix}W_i =\\ \ \end{matrix} \ \begin{matrix}\text{the number of $i$-dimensional linear subspaces $V'$ in $V$ such}\\ \text{that $V'$ is the span of a subset of the vectors $\{v_1, \ldots, v_n\}$.\quad}\end{matrix}
\]
\end{enumerate}

We leave it as an exercise to check that for the matroid associated to the graph or the vector configuration in \Cref{fig:intro}, we have:

\smallskip
\begin{tabular}{ll}
$(I_0, I_1, I_2, I_3) = (1,5,10,8)$, \qquad &$(w_0, w_1, w_2, w_3) = (1,5,8,4)$,\\
$(I'_0, I'_1, I'_2, I'_3) = (1,2,3,2)$, \qquad &$(w'_0, w'_1, w'_2) = (1,2,1)$,\\
$(W_0, W_1, W_2, W_3) = (1,5,6,1)$.
\end{tabular}

\medskip
\noindent Notice in this example that every sequence is unimodal, log-concave, and top-heavy.
Several conjectures from the 70's posited that these sequences are unimodal, log-concave, or top-heavy for an arbitrary matroid.  We describe these conjectures and their history more fully in \Cref{subsec:Tutte}.

\subsection{An approach from algebraic geometry}
After decades of little progress, a breakthrough occurred when many of these conjectures were resolved for realizable matroids using algebraic geometry: Huh and Katz \cite{Huh12, HK12} showed that the sequence \ref{seqw} is log-concave (with no internal zeros), Huh \cite{Huh15} showed that \ref{seqw'} is log-concave (with no internal zeros), and Huh and Wang \cite{HW17} showed that \ref{seqW} is top-heavy.
These developments were particularly significant in light of the following phenomena in matroid theory:

\smallskip
The geometry of realizable matroids often inspires purely combinatorial constructions for all matroids.
Certain geometric properties, a priori applicable only to realizable matroids, persist to all matroids through these purely combinatorial constructions.
This is surprising because almost all matroids are not realizable \cite{Nel18}, but 
such a creative tension between geometry and combinatorics is a recurring theme in matroid theory.

\smallskip
A recent spectacular example of this phenomenon is the development of the Hodge theory of matroids by June Huh and his collaborators \cite{AHK18, ADH22, BHMPW22, BHMPWb}, which successfully resolved conjectures about log-concavity or top-heaviness of the sequences \ref{seqI},\ldots,\ref{seqW} for arbitrary (not necessarily realizable) matroids.
They established that matroids satisfy combinatorial analogues of certain Hodge-theoretic properties in algebraic geometry, known sometimes as the ``K\"ahler package'':

\begin{defn}\label{defn:Kahler}
Let $A^\bullet = \bigoplus_{i=0}^d A^i$ be a finite-dimensional graded real vector space with a symmetric bilinear form $P: A^\bullet \times A^{d-\bullet} \to \RR$, and let $\mathcal K$ be a convex subset  of graded linear operators $L: A^\bullet \to A^{\bullet+1}$ satisfying $P(Lx,y) = P(x,Ly)$ for all $x,y\in A^\bullet$.  The triple $(A^\bullet, P, \mathcal K)$ is said to satisfy the \textbf{K\"ahler package} if the following three properties hold for all nonnegative integers $i\leq \frac{d}{2}$:
\begin{itemize}
\item[(PD)] \label{PD} The pairing $P: A^i \times A^{d-i} \to \RR$ is non-degenerate (Poincar\'e duality).
\item[(HL)] \label{HL} For any $L_1, \ldots, L_{d-2i} \in \mathcal K$, the linear map
\[
A^{i} \to A^{d-i} \quad\text{given by}\quad x\mapsto L_1\cdots L_{d-2i}x\]
is an isomorphism (hard Lefschetz property in degree $i$).
\item[(HR)] \label{HR} For any $L_0,L_1, \ldots, L_{d-2i} \in \mathcal K$, the symmetric bilinear pairing
\[
A^i \times A^i \to \RR \quad\text{given by}\quad (x,y) \mapsto (-1)^iP(x,L_1\cdots L_{d-2i}y)
\]
is positive definite when restricted to the kernel of the map $A^i\to A^{d-i+1}$ given by $x\mapsto L_0L_1\cdots L_{d-2i} x$ (Hodge-Riemann relations in degree $i$).
\end{itemize}
\end{defn}

Classical Hodge theory tells us that these properties are satisfied when $A^\bullet$ is the cohomology ring of real $(p,p)$-forms on a complex projective manifold, $P$ is the Poincar\'e duality pairing, and $\mathcal K$ consists of multiplication by ample divisor classes (see \cite{Huy05} or \cite{Voi02}).

\smallskip
The geometry behind realizable matroids led to purely combinatorial constructions for various ``cohomologies'' of a matroid.  These constructions include the Chow ring of a matroid \cite{FY04, AHK18}, the conormal Chow ring of a matroid \cite{ADH22}, and the intersection cohomology of a matroid \cite{BHMPWb}.   For a matroid realizable over $\CC$, all three satisfy the K\"ahler package due to classical algebraic geometry.  The incredible result of June Huh and his collaborators---Karim Adiprasito, Federico Ardila, Tom Braden, Graham Denham, Eric Katz,  Jacob Matherne, Nick Proudfoot, Botong Wang---is that the K\"ahler package continues to hold for these ``cohomologies'' of a matroid even when the matroid is not realizable.  

\medskip
We survey this remarkable development in matroid theory and its connection to algebraic geometry in four parts.  In \Cref{sec:matroids}, we give a brief introduction to matroids, and describe the long-standing conjectures resolved by the Hodge theory of matroids.  In \Cref{sec:real}, we explain how the conjectures in the case of realizable matroids were resolved using algebraic geometry.
In \Cref{sec:tropHodge}, we discuss the K\"ahler package for Chow rings of fans and matroids, and how the validity of (HR) implies the conjectures about log-concavity.  In \Cref{sec:IH}, we discuss the intersection cohomology of a matroid, and explain its implication to top-heaviness.

\medskip
Several interesting topics had to be omitted, even though they are closely related to the topics discussed here.  A partial list includes the following:
\begin{itemize}
\item The Kazhdan-Lusztig theory of matroids \cite{EPW16}, which was an inspiration behind the construction of intersection cohomology of a matroid.  We point to \cite{Pro18} for a survey of Kazhdan-Lusztig-Stanley polynomials in a more general context.
\item The study of matroids through the polyhedral geometry of their base polytopes, their subdivisions, and the geometry of the Grassmannian \cite{GGMS87, Laf03}.  We point to the survey \cite[Section 4]{Ard22} and references therein.

\item Other approaches to the ``cohomology'' of a matroid in the broader context of tropical geometry, for instance \cite{IKMZ19, AP}.
\end{itemize}
We hope that this survey will spark the reader's general interest in this active field of the study of matroids from an algebro-geometric perspective.

\proof[Notation]
Throughout, let $E = \{1, \ldots, n\}$ be a finite set of cardinality $n$.  For a subset $S\subseteq E$, we denote by $\be_S = \sum_{i \in S} \be_i$ the sum of standard basis vectors in $\kk^E$, where the field $\kk$ will be clear in context.
An algebraic variety is reduced and irreducible (over an algebraically closed field).

\subsection*{Acknowledgements} {The author thanks June Huh for helpful conversations, and Matt Baker, Andrew Berget, Tom Braden, Eric Katz, Matt Larson, and Lauren Williams for helpful suggestions and comments on a preliminary draft of this article.}

\section{Background in matroid theory}\label{sec:matroids}

Here we give a minimal introduction to matroids.
In addition to standard references on matroids such as \cite{Wel76, Oxl11}, we point to \cite{Ard22, Bak18, Huh18b, Kat16} for surveys tailored towards studying matroids from an algebro-geometric viewpoint.

\subsection{Constructions}\label{subsec:basic}

Since subsets of independent sets are independent, we may specify a matroid by its maximal independent sets, called the \textbf{bases} of the matroid.

\begin{eg}
For an integer $0\leq r \leq n$, the \textbf{uniform matroid} of rank $r$ on $E$ is the matroid $\mathrm{U}_{r,n}$ whose bases are all subsets of cardinality $r$.  When $n=r$, we say that $\mathrm{U}_{n,n}$ is the \textbf{Boolean matroid} on $E$.  The Boolean matroid $\U_{0,0}$, i.e.\ when $E = \varnothing$ so $n=r=0$, is called the \textbf{trivial matroid}.
Any uniform matroid $\mathrm{U}_{r,n}$ is realizable over any infinite field $\kk$ as a general collection of $n$ vectors in $V = \kk^r$.
\end{eg}

\begin{eg}\label{eg:cover}
We may visualize a collection of vectors in a 4-dimensional vector space as a collection of points in the projective 3-space $\PP^3$.  For example, the 5 column vectors of the matrix

\begin{figure}[h]
	\centering
\begin{tabular}{cc}

\raisebox{37pt}{$\begin{bmatrix}
1 & 0 & 0 & -1 & 0 \\
0 & 1 & 0 & -1 & 0 \\
0 & 0 & 1 & -1 & 0 \\
0 & 0 & 0 & 0 & 1 \\
\end{bmatrix}$}

&
\qquad\quad\includegraphics[height=27mm]{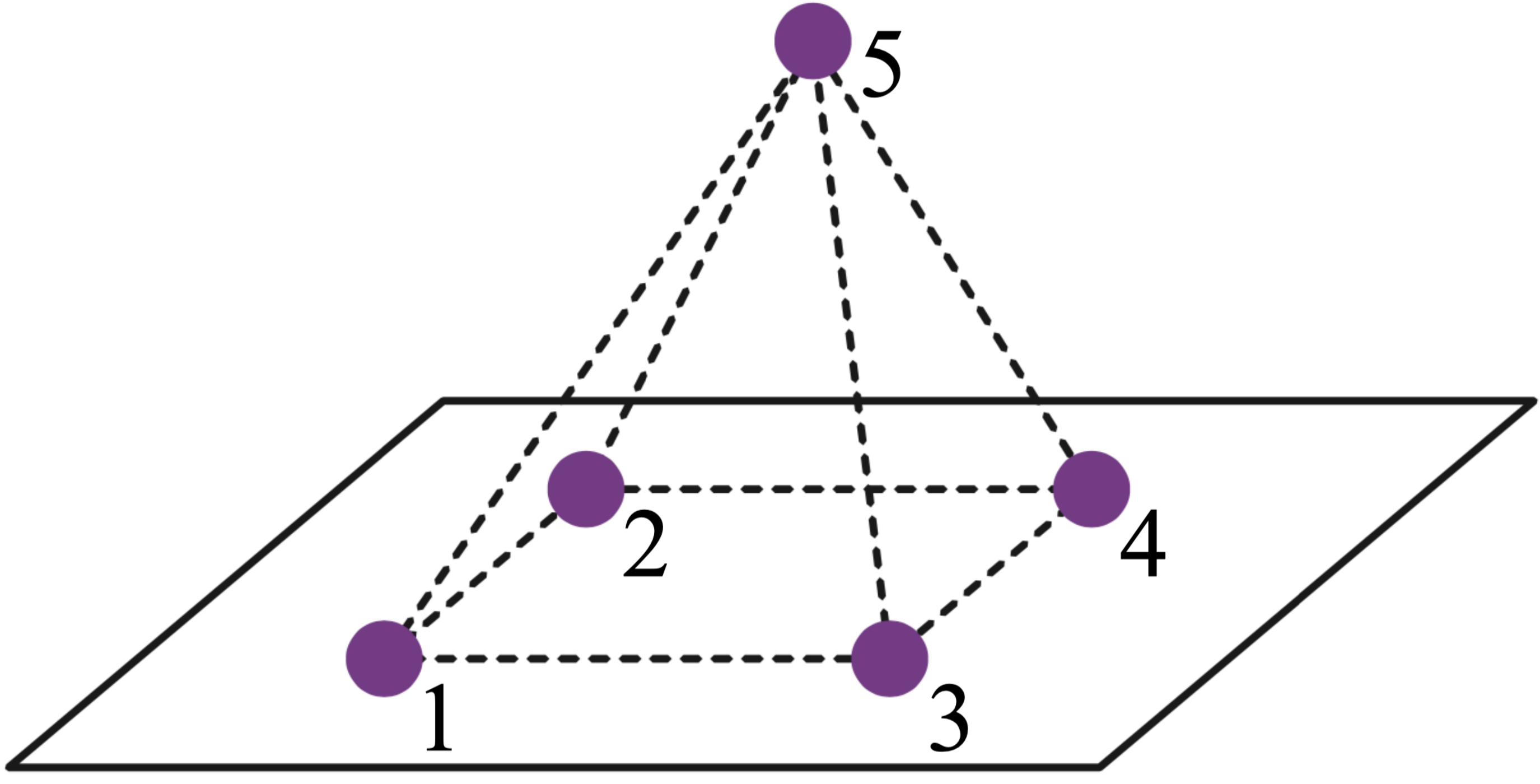}
\vspace{-15pt}

\end{tabular}
\end{figure}

\noindent can be visualized as the purple points, four of which lie in a common projective plane.  The bases of this matroid are $\{1235,1245,1345,2345\}$.
\end{eg}

We define the \textbf{dual matroid} $\M^\perp$ of a matroid $\M$ on ground set $E$ by declaring
\[
\text{the set of bases of $\M^\perp$} = \{E\setminus B : B \text{ a basis of } \M\}.
\]
For example, we have $\U_{r,n}^\perp = \U_{n-r,n}$.  For the matroid $\M$ in \Cref{eg:cover}, its dual $\M^\perp$ has the set of bases $\{1,2,3,4\}$.
Many notions in matroid theory come in pairs via matroid duality.  For instance, an element $e\in E$ is a \textbf{loop} of $\M$ if it is in no bases, and is a \textbf{coloop} if it is in every basis of $\M$.
The matroid in \Cref{eg:cover} has no loops and has a coloop 5, or equivalently, its dual matroid has no coloops and has a loop 5.

\medskip
Another useful way of describing a matroid is by its rank function.
For a matroid $\M = (E,\mathcal  I)$, its \textbf{rank function} $\rk_\M: 2^E \to \ZZ$ is defined by
\[
\rk_\M(S) = \max\{ |I| : I\in \mathcal  I \text{ and } I\subseteq S\} \quad\text{for every subset $S\subseteq E$}.
\]
In particular, an independent set of a matroid $\M$ is a subset $I\subseteq E$ whose rank $\rk_\M(I)$ equals its cardinality $|I|$.
That is, independent sets are the minimal subsets of $E$ with respect to a given rank.  Considering the maximal subsets leads to the notion of flats of a matroid.

\begin{defn}
A \textbf{flat} of a matroid $\M$ on $E$ is a subset of $E$ that is maximal for its rank.  That is, a subset $F\subseteq E$ is a flat of $\M$ if
$\rk_\M(F\cup \{e\}) > \rk_\M(F) \text{ for all $e\in E\setminus F$}$.
\end{defn}

The set of flats of a matroid $\M$ forms a poset under inclusion.  This poset is a \emph{lattice} with meet and join defined by
\[
F \wedge F' = F\cap F' \quad\text{and}\quad F\vee F' = \text{the smallest flat containing $F\cup F'$}.
\]
For example, the lattice of flats of the matroid in \Cref{eg:cover} is

\begin{center}
\begin{tikzpicture}[scale=0.85, vertices/.style={draw, fill=black, circle, inner sep=0pt}]
              \node [vertices, label=right:{$\varnothing$}] (0) at (-0+0,0){};
              \node [vertices, label=right:{$5$}] (1) at (-3+0,4/3){};
              \node [vertices, label=right:{$4$}] (2) at (-3+3/2,4/3){};
              \node [vertices, label=right:{$3$}] (3) at (-3+3,4/3){};
              \node [vertices, label=right:{$2$}] (4) at (-3+9/2,4/3){};
              \node [vertices, label=right:{$1$}] (5) at (-3+6,4/3){};
              \node [vertices, label=left:{$45$}] (6) at (-27/4+0,8/3){};
              \node [vertices, label=left:{$34$}] (8) at (-27/4+3/2,8/3){};
              \node [vertices, label=left:{$35$}] (7) at (-27/4+3,8/3){};
              \node [vertices, label=left:{$25$}] (9) at (-27/4+9/2,8/3){};
              \node [vertices, label=left:{$24$}] (10) at (-27/4+6,8/3){};
              \node [vertices, label=left:{$23$}] (13) at (-27/4+15/2,8/3){};
              \node [vertices, label=right:{$15$}] (11) at (-27/4+9,8/3){};
              \node [vertices, label=right:{$24$}] (12) at (-27/4+21/2,8/3){};
              \node [vertices, label=right:{$13$}] (14) at (-27/4+12,8/3){};
              \node [vertices, label=right:{$12$}] (15) at (-27/4+27/2,8/3){};
              \node [vertices, label=left:{$345$}] (16) at (-9/2+0,4){};
              \node [vertices, label=right:{$245$}] (17) at (-9/2+3/2,4){};
              \node [vertices, label=right:{$235$}] (19) at (-9/2+3,4){};
              \node [vertices, label=right:{$145$}] (18) at (-9/2+9/2,4){};
              \node [vertices, label=right:{$135$}] (20) at (-9/2+6,4){};
              \node [vertices, label=right:{$125$}] (21) at (-9/2+15/2,4){};
              \node [vertices, label=right:{$1234$}] (22) at (-9/2+9,4){};
              \node [vertices, label=right:{$12345$}] (23) at (-0+0,16/3){};
      \foreach \to/\from in {0/1, 0/2, 0/3, 0/4, 0/5, 1/9, 1/11, 1/6, 1/7, 2/8, 2/10, 2/12, 2/6, 3/8, 3/13, 3/14, 3/7, 4/9, 4/10, 4/13, 4/15, 5/11, 5/12, 5/14, 5/15, 6/16, 6/17, 6/18, 7/20, 7/16, 7/19, 8/16, 8/22, 9/21, 9/17, 9/19, 10/17, 10/22, 11/20, 11/21, 11/18, 12/18, 12/22, 13/22, 13/19, 14/20, 14/22, 15/21, 15/22, 16/23, 17/23, 18/23, 19/23, 20/23, 21/23, 22/23}
      \draw [-] (\to)--(\from);
      \end{tikzpicture}
\end{center}

\begin{exer}\
\begin{enumerate}
\item Show that every subset is contained in a unique flat of the same rank.  In particular, the join of flats is well-defined.
\item How can one recover the bases of a matroid from the lattice of its flats?
\item Compute the lattice of flats of the matroid of the graph in \Cref{fig:intro}.
\end{enumerate}
\end{exer}

We record some linear algebraic interpretations of the notions introduced for a matroid $\M$ realized by a set of vectors $\{v_1, \ldots, v_n\}$ spanning a vector space $V$:
\begin{itemize}
\item The bases of $\M$ are subsets $B\subseteq E$ such that $\{v_i : i\in B\}$ is a basis of $V$.
\item An element $e\in E$ is a loop if and only if $v_e =0$.
\item We have $\rk_\M(S) = \dim \operatorname{span}\{v_i : i \in S\}$ for any subset $S\subseteq E$.
\item The flats of $\M$ are subsets $F\subseteq E$ such that $F = V' \cap \{v_1, \ldots, v_n\}$ for some linear subspace $V'\subseteq V$.  That is, the flats correspond to the different spans of subsets of the vectors.
\end{itemize}

Per the last bullet point, we now define the sequence \ref{seqW} for arbitrary matroids.

\begin{defn}
For a matroid $\M$ of rank $r$, the \textbf{Whitney numbers of the second kind} $(W_0, W_1, \ldots, W_r)$ are defined by
\[
W_i = \text{the number of flats of $\M$ with rank $i$}.
\]
\end{defn}

\begin{exer}
Show that if $\M$ is the graphical matroid of the complete graph on $N$ vertices, then its flats of rank $i$ correspond to partitions of $N$ into $N-i$ (nonempty) parts.  In particular, the numbers $(W_0, W_1, \ldots, W_{N-1})$ in this case are known as the Stirling numbers of the second kind.
\end{exer}

Matroid duality also admits a linear algebraic interpretation, given in the next exercise.

\begin{exer}\label{exer:dual}
Let $\kk^E \twoheadrightarrow V$ be the map given by $\be_i \mapsto v_i$, and let $K$ be its kernel.  The short exact sequence $0\to K \to \kk^E \to V \to 0$ dualizes to $0\to V^\vee \to \kk^E \to K^\vee \to 0$.  Show that the surjection $\kk^E \twoheadrightarrow K^\vee$ realizes the dual matroid $\M^\perp$.
\end{exer}

An important pair of operations for forming new matroids from a given matroid is the restriction and contraction.

\begin{defn}
For a matroid $\M$ on ground set $E$, and a subset $A\subseteq E$, we define two matroids $\M|A$ and $\M/A$ on ground sets $A$ and $E\setminus A$, respectively, by specifying their rank functions:
\begin{align*}
&\rk_{\M|A}(S) = \rk_\M(S) &\text{for all $S\subseteq A$.\qquad}\\
&\rk_{\M/A}(S) = \rk_\M(S\cup A) - \rk_\M(A) &\text{for all $S\subseteq E\setminus A$}.
\end{align*}
The matroid $\M|A$ is called the \textbf{restriction} of $\M$ to $A$, and the matroid $\M/A$ is called the \textbf{contraction} of $\M$ by $A$.  The \textbf{deletion} $\M\setminus A$ of $\M$ by $A$ is the restriction $\M|({E\setminus A})$.
\end{defn}

These operations behave particularly well for a flat $F$ of a matroid $\M$:
\begin{itemize}
\item The set of flats of $\M|F$ is $\{F' : F' \text{ a flat of $\M$ contained in $F$}\}$.
\item The set of flats of $\M/F$ is $\{F'\setminus F : F' \text{ a flat of }\M \text{ containing $F$}\}$.
\end{itemize}

These operations have a graphical and linear algebraic interpretations as well:
For the graphical matroid of a graph $G$, deletion corresponds to deleting the corresponding edges, and contraction corresponds to contracting the corresponding edges.
When a matroid $\M$ is realized by a set of vectors $\{v_i : i\in E\}$, the restriction $\M|A$ is realized by the subset of vectors $\{v_i: i\in A\}$.  The contraction $\M/A$ is realized by the images of the vectors $\{v_i: i \in E\setminus A\}$ under the quotient by the span of $\{v_i: i\in A\}$.

\begin{exer}
Show that deletion and contraction are dual notions, that is, we have $(\M\setminus A)^\perp = \M^\perp / A$.
\end{exer}

\subsection{Invariants}\label{subsec:Tutte}

Introduced for graphs by Tutte \cite{Tut67} and extended to matroids by Crapo \cite{Cra69}, the Tutte polynomial is among the most famous invariants of a matroid.  For the proofs of the statements here, as well as a fuller treatment of Tutte polynomials, see \cite{BO92}.

\begin{defn}
Tutte polynomial of a matroid $\M$ of rank $r$ on ground set $E$ is the bivariate polynomial defined by
\[
\T_\M(x,y)  = \sum_{S\subseteq E} (x-1)^{r - \rk_\M(S)}(y-1)^{|S|-\rk_\M(S)}.
\]
\end{defn}

The Tutte polynomial is the universal deletion-contraction invariant:

\begin{thm}\label{thm:delcont}
The Tutte polynomial can be defined recursively by
\[
\T_\M(x,y) = \begin{cases}
x \T_{\M/ e}(x,y) & \text{if $e\in E$ a coloop in $\M$}\\
y \T_{\M\setminus e}(x,y) & \text{if $e\in E$ a loop in $\M$}\\
\T_{\M\setminus e}(x,y) + \T_{\M/e}(x,y) &\text{if $e\in E$ neither loop nor coloop}
\end{cases}
\]
with $\T_{\U_{0,0}}(x,y) = 1$.  If $f$ is an invariant of matroids with values in a (commutative unital) ring $R$ such that $f(\U_{0,0})= 1$ and there exists $x_0, y_0, a, b\in R$ satisfying
\[
f(\M) = \begin{cases}
x_0 f({\M/ e}) & \text{if $e\in E$ a coloop in $\M$}\\
y_0 f(\M\setminus e) & \text{if $e\in E$ a loop in $\M$}\\
af({\M\setminus e}) + bf({\M/e}) &\text{if $e\in E$ neither loop nor coloop}
\end{cases}
\]
for all matroids $\M$ and an element $e$.  Then, we have
\[
f(\M) = a^{|E|-\rk_\M(E)}b^{\rk_\M(E)} \T_\M\left(\frac{x_0}{b}, \frac{y_0}{a}\right).
\]
\end{thm}

The theorem implies the following basic properties of the Tutte polynomial:
\begin{itemize}
\item The Tutte polynomial $T_\M(x,y)$ of a matroid $\M$ has nonnegative coefficients.
\item The constant term of $T_\M(x,y)$ is zero unless $\M$ is the trivial matroid.
\item For the dual matroid $\M^\perp$, we have $T_{\M^\perp}(x,y) = T_\M(y,x)$.
\end{itemize}

\medskip
Univariate specializations of the Tutte polynomial lead to many interesting combinatorial sequences.  For example, we deduce from the definition that for a matroid $\M$ of rank $r$,
\[
T_\M(q+1,1) = \sum_{i=0}^r I_i q^{r-i},
\]
where we recall from \ref{seqI} that $I_i$ is the number of independent sets of $\M$ of cardinality $i$.  Consequently, we have that the sequence $(I_0', \ldots, I_r')$ of \ref{seqI'} is given by
\[
T_\M(q,1) = \sum_{i=0}^r I_i' q^{r-i}.
\]

Another important specialization is the \textbf{characteristic polynomial} $\chi_\M$ of a matroid $\M$ of rank $r$, defined as
\[
\chi_\M(q) = (-1)^r T_\M(1-q,0).
\]
The basic properties of the Tutte polynomial listed above imply that the coefficients of $\chi_\M(q)$ have alternating signs, and that $\chi_\M(q)$ is divisible by $(q-1)$ unless $\M$ is a trivial matroid.
Thus, one often divides out $(q-1)$ to define the \textbf{reduced characteristic polynomial} $\overline\chi_\M(q) = \chi_\M(q)/(q-1)$.
It follows from \Cref{thm:delcont} that $\chi_\M = 0$ if $\M$ has a loop.

\medskip
The characteristic polynomial of a graphical matroid essentially equals the chromatic polynomial of the graph in the following way.

\begin{exer}
Show that if $\M$ is the graphical matroid of a finite graph $G$, and $G$ has $c$ connected components, then $q^c \chi_\M(q) = \chi_G(q)$.  (Hint: Appeal to \Cref{thm:delcont} by showing that both satisfy an identical deletion-contraction relation).
\end{exer}

We may now define the sequences \ref{seqw} and \ref{seqw'}, which we previously only defined for graphical matroids, for arbitrary nontrivial matroids:  they are the coefficients of $T_\M(1+q,0)$ and $\frac{1}{q}T_\M(q,0)$, respectively.

\medskip
Summarizing, we have the following sequences for a matroid $\M$ of rank $r$.
\begin{enumerate}[label=(\alph*)]
\item \label{sseqI} The coefficients $(I_0, \ldots, I_r)$ of $T_\M(q+1,1) = \sum_{i=0}^r I_i q^{r-i}$.
\item \label{sseqI'} The coefficients $(I'_0, \ldots, I'_r)$ of $T_\M(q,1) = \sum_{i=0}^r I'_i q^{r-i}$.
\item \label{sseqw}The coefficients $(w_0, \ldots, w_r)$ of $T_\M(q+1,0) = \sum_{i=0}^r w_i q^{r-i}$.
\item \label{sseqw'} The coefficients $(w'_0, \ldots, w'_{r-1})$ of $\frac{1}{q}T_\M(q,0) = \sum_{i=0}^r I_i q^{r-1-i}$.
\item \label{sseqW} The Whitney numbers of the second kind $(W_0, \ldots, W_r)$ of $\M$.
\end{enumerate}

\begin{thm}\cite{AHK18, ADH22, BHMPWb} \label{thm:main1} 
Let $\M$ be a matroid of rank $r$.
\begin{enumerate}[label = (\alph*)]
\item The sequence $(I_0, \ldots, I_r)$ is unimodal, log-concave, and top-heavy.
\item The sequence $(I'_0, \ldots, I'_r)$ is unimodal, log-concave, and top-heavy.
\item The sequence $(w_0, \ldots, w_r)$ is unimodal, log-concave, and top-heavy.
\item The sequence $(w'_0, \ldots, w'_{r-1})$ is unimodal, log-concave, and top-heavy.
\item The sequence $(W_0, \ldots, W_r)$ satisfies $W_i \leq W_j$ for all $0\leq i \leq j \leq r -i$.  In particular, it is top-heavy.
\end{enumerate}
\end{thm}

The statements of the theorem were long-standing conjectures in matroid theory.  The unimodality and log-concavity conjectures are due to: Welsh \cite{Wel71} and Mason \cite{Mas72} for \ref{sseqI}, Dawson \cite{Daw84} for \ref{sseqI'}, Read \cite{Rea68} and Hoggar \cite{Hog74} for \ref{sseqw} of graphical matroids, Heron \cite{Her72}, Rota \cite{Rot71}, and Welsh \cite{Wel76} for \ref{sseqw}, and Brylawski \cite{Bry82} for \ref{sseqw'}.
Hibi \cite{Hib92} and Swartz \cite{Swa03} posed the top-heaviness of \ref{sseqI'} and \ref{sseqw'} (respectively).\footnote{A slightly different terminology of \emph{flawless-ness} appears in \cite{Hib92} and related works \cite{Hib89, JKL18}.  A nonnegative sequence $(a_0, \ldots, a_m)$ is flawless if it is top-heavy and additionally satisfies $a_0 \leq \cdots \leq a_{\lfloor {d}/{2}\rfloor}$, where $d$ is the largest index such that $a_d \neq 0$.  Note that unimodal and top-heavy sequences are flawless.}
Dowling and Wilson \cite{DW74, DW75} conjectured the top-heaviness of \ref{sseqW}, generalizing a theorem of de Bruijn and Erd\"os \cite{dBE48} on point-line incidences in projective planes.
There are two notable conjectures on \ref{sseqW} that remain open: Rota \cite{Rot71} conjectured its unimodality, and Mason \cite{Mas72} its log-concavity.

\begin{rem}
One may ask if there is a log-concavity statement for the whole Tutte polynomial of a matroid that explains the log-concavity of the four specializations in \Cref{thm:main1}.  This was achieved by Berget, Spink, Tseng, and the author \cite{BEST} who showed that the 4-variable transformation
\[
(x+y)^{-1}(y+z)^r(x+w)^{|E|-r}\T_\M\left(\frac{x+y}{y+z}, \frac{x+y}{x+w}\right)
\]
of the Tutte polynomial of a matroid $\M$ satisfies a multivariate version of log-concavity.
We note that without such a transformation, the coefficients of $\T_\M(x,y)$ can fail to be unimodal \cite{Sch93}.
\end{rem}

In this survey, we will explain how the sequences \ref{sseqw} and \ref{sseqw'} are shown to be log-concave with no internal zeros,\footnote{As observed by Lenz \cite{Len13}, a result of Brylawski \cite[Theorem 4.2]{Bry77} implies that the statements for \ref{sseqI} and \ref{sseqI'} follow from those for \ref{sseqw} and \ref{sseqw'}.  The top-heaviness of \ref{sseqw} and \ref{sseqw'} follows from their unimodality due to \cite[Theorem 1.2]{JKL18}.
}
and how the sequence \ref{sseqW} is shown to be top-heavy.  
Since $\chi_\M = 0$ if $\M$ has a loop, and deleting loops of a matroid does not change the lattice of its flats, we assume the following:

\medskip
\noindent \textbf{Assumption.}  From now on, a matroid is loopless unless specified otherwise.

\section{The realizable case}\label{sec:real}

We explain how the statements in \Cref{thm:main1} can be deduced using algebraic geometry when the matroid in question is realizable.
We assume familiarity with algebraic geometry; those who prefer purely combinatorial treatments may skip this section.
For simplicity, we consider matroids realizable over $\CC$.  For other fields, one can run nearly identical arguments using the Chow cohomology ring \cite{Ful98} or the $\ell$-adic (intersection) cohomology in place of singular (intersection) cohomology $(I)H^\bullet(-)$ with rational coefficients.

\medskip
Throughout this section, let $\M$ be a (nontrivial) matroid of rank $r$ realized by a set of vectors $\{v_i : i\in E\}$ spanning a vector space $V \simeq \CC^r$.  The corresponding surjection $\CC^E \twoheadrightarrow V$ dualizes to give an $r$-dimensional linear subspace $V^\vee\subseteq \CC^E$.

\medskip
\noindent \textbf{Notation.} Let us denote $L = V^\vee$ to avoid repeated use of the superscript $\vee$.

\medskip
The set of independent sets of $\M$ can then be described also as the collection
\[
\mathcal I = \{I\subseteq E: \text{the composition $L\hookrightarrow \kk^E \twoheadrightarrow \kk^I$ is surjective}\}.
\]
We will often projectivize and work with $\PP L \subseteq \PP^{n-1}$.

\subsection{Hyperplane arrangements}
We first discuss some structures of the matroid $\M$ in terms of its realization as a subspace $L\subseteq \CC^E$.
For each $i\in E$ let $H_i$ be the $i$-th coordinate hyperplane of $\CC^E$.  Our assumption that $\M$ is loopless implies that $L$ is not contained in any $H_i$.  We thus have a \textbf{hyperplane arrangement} $\mathcal A$ on $L$ consisting of the hyperplanes $\{L\cap H_i: i\in E\}$.  Dualizing the correspondence between the flats of $\M$ and the spans of subsets of the vectors $\{v_i: i\in E\}$, we obtain a correspondence
\begin{align*}
\{\text{flats of $\M$}\} &\longleftrightarrow \{\text{subspaces of $L$ arising as intersections of hyperplanes in $\mathcal A$}\}\\
F &\longleftrightarrow L_F = L\cap \bigcap_{i\in F}H_i.
\end{align*}
Note that the correspondence is order-reversing, and in particular, a flat of rank $r-i$ maps to the linear subspace $L_F$ of dimension $i$.

\begin{eg}\label{eg:U34}
The columns of the matrix below realizes the matroid $\U_{3,4}$. 
Equivalently, the embedded subspace $L\subseteq \CC^4$, where $L = \{x_1 + x_2 + x_3 + x_4 = 0 \}$ is the row-span of the matrix, realizes the matroid $\U_{3,4}$.

\begin{figure}[h]
	\centering
\begin{tabular}{ccc}

\raisebox{50pt}{
$\begin{bmatrix}
1 & 0 & 0 & -1 \\
0 & 1 & 0 & -1 \\
0 & 0 & 1 & -1 \\
\end{bmatrix}
$}

&
\raisebox{10pt}{
\qquad \begin{tikzpicture}[scale=0.65, vertices/.style={draw, fill=black, circle, inner sep=0pt}]
              \node [vertices] (0) at (-0+0,0){};
              \node [vertices] (1) at (-2.25+0,1.33333){};
              \node [vertices] (2) at (-2.25+1.5,1.33333){};
              \node [vertices] (3) at (-2.25+3,1.33333){};
              \node [vertices] (4) at (-2.25+4.5,1.33333){};
              \node [vertices] (5) at (-3.75+0,2.66667){};
              \node [vertices] (8) at (-3.75+1.5,2.66667){};
              \node [vertices] (6) at (-3.75+3,2.66667){};
              \node [vertices] (9) at (-3.75+4.5,2.66667){};
              \node [vertices] (10) at (-3.75+6,2.66667){};
              \node [vertices] (7) at (-3.75+7.5,2.66667){};
              \node [vertices] (11) at (-0+0,4){};
      \foreach \to/\from in {0/1, 0/2, 0/3, 0/4, 1/5, 1/6, 1/7, 2/8, 2/5, 2/9, 3/8, 3/6, 3/10, 4/9, 4/10, 4/7, 5/11, 6/11, 7/11, 8/11, 9/11, 10/11}
      \draw [-] (\to)--(\from);
      \end{tikzpicture}
}

&
\includegraphics[height=32mm]{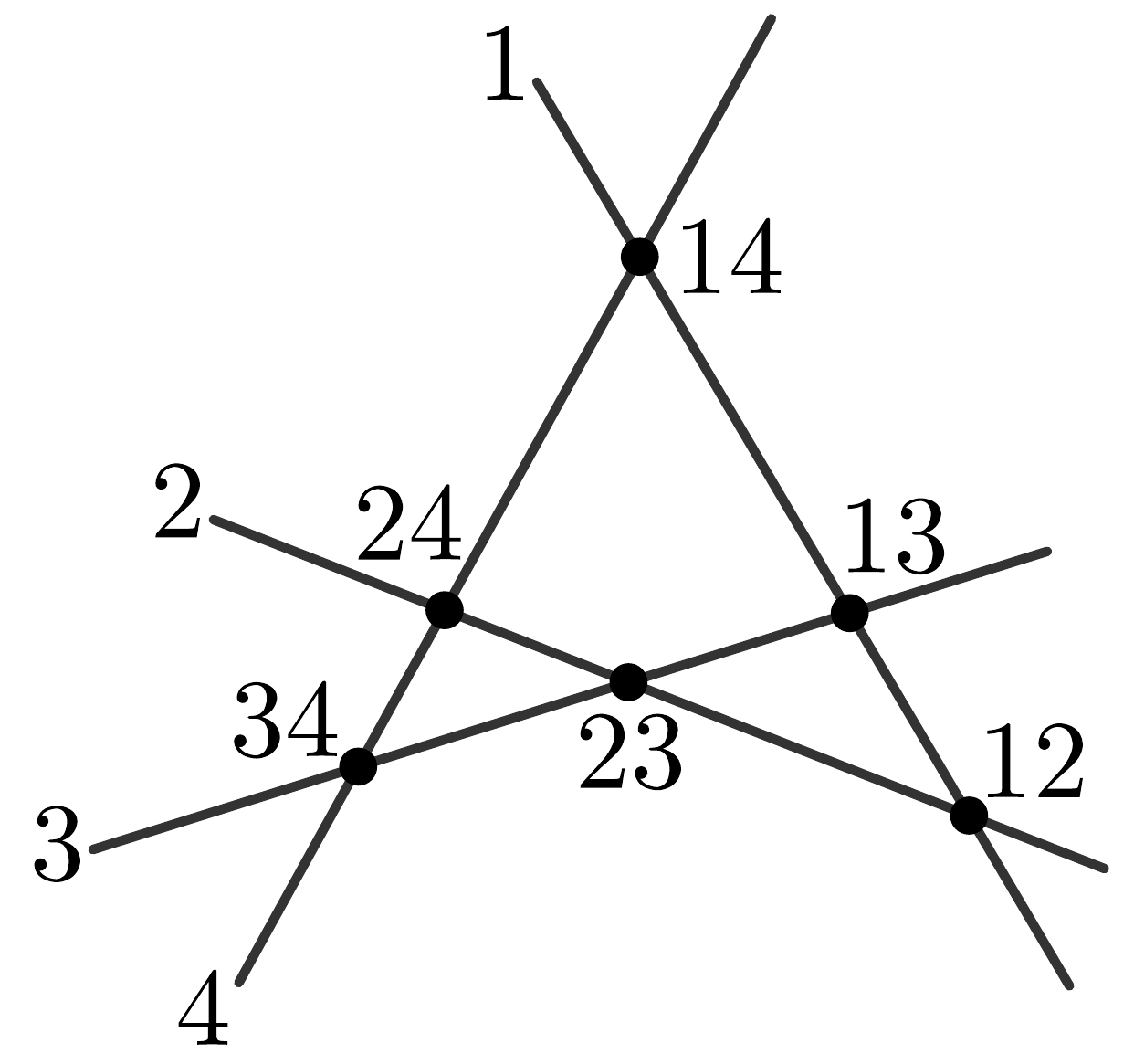}
\end{tabular}
\end{figure}
\noindent Next to the matrix, we have depicted the lattice of flats of $\M$ and the projectivization of the hyperplane arrangement $\mathcal A$ in $\PP L \simeq \PP^2$.
\end{eg}

We denote the complement of the hyperplane arrangement by
\[
\mathring L =  \textstyle L \setminus \bigcup \mathcal A = L \cap (\CC^*)^E, \quad\text{and likewise,}\quad \PP \mathring L = \PP L \cap \big((\CC^*)^E/\CC^*\big).
\]

\begin{exer}\label{exer:scissors}
Suppose $i\in E$ is not a coloop, and let $F$ be the smallest flat containing $i$.  Show that the subspace $L_F$ is a realization of the contraction $\M/F$, and that the hyperplane arrangement complement $ L \setminus \bigcup(\mathcal A \setminus \{L_F\})$ is a realization of the deletion $\M\setminus F$.  What happens when $i$ is a coloop?
\end{exer}

The geometric study of (complements of) hyperplane arrangements and its interaction with matroid theory is a rich and on-going research field; some references include \cite{OT92, Dim17}.  Here, we only note the following classical fact \cite{OS80}, which states that the characteristic polynomial records the dimensions of the cohomologies of the arrangement complement.

\begin{thm}\label{thm:poin}
Let $L\subseteq \CC^E$ realize a matroid $\M$ of rank $r$.  Then, we have
\[
\T_\M(1+q,0) = \sum_{i=0}^r \dim H^i(\mathring L) q^{r-i},
\]
or equivalently, by applying the K\"unnuth formula to $\PP \mathring L \times \CC^* \simeq \mathring L$,
\[
\frac{1}{1+q}\T_\M(1+q,0) = \sum_{i=0}^{r-1} \dim H^i(\PP \mathring L) q^{r-1-i}.
\]
\end{thm}

\subsection{Log-concavity via intersection degrees}
We now explain how the log-concavity of the sequences \ref{sseqw} and \ref{sseqw'} can be shown using algebraic geometry.
We start by describing a general strategy for log-concavity.

\medskip
Let $X$ be a smooth projective $\CC$-variety $X$ of dimension $d$.  When considered as a $2d$-dimensional real compact manifold, Poincar\'e duality for  the cohomology ring $H^\bullet(X)$ provides the isomorphism $\int_X: H^{2d}(X) \to \ZZ$, called the degree map.
Recall that a divisor $D$ on $X$ is ample (resp.\ semi-ample) if the line bundle $\mathcal O_X(mD)$ is very ample (resp.\ globally generated) for some integer $m>>0$.
The general strategy arises from the following Khovanskii-Teissier inequalities (see \cite[Section 1.6]{Laz04a} for a history and a fuller treatment).

\begin{prop}\label{prop:KT}
Let $\alpha,\beta \in H^2(X)$ be the cohomology classes of two semi-ample (or more generally, nef) divisors on a smooth projective $\CC$-variety $X$ of dimension $d$. Then,
\[\text{the sequence $(a_0, \ldots, a_d)$ of intersection degrees of $\alpha$ and $\beta$, i.e.\ } a_i = \int_X \alpha^{d-i}\beta^i
\]
is log-concave with no internal zeros.
\end{prop}

\begin{proof}[Sketch of the proof]
By continuity, one can assume $\alpha,\beta$ to be ample, and then one reduces to the case of surfaces via the Bertini theorem.  Then, one of the equivalent forms of the Hodge index theorem for surfaces \cite[Exercise V.1.9]{Har77} exactly yields the desired log-concavity.  We note that the Hodge index theorem for surfaces stated as \cite[Theorem V.1.9]{Har77} is exactly the validity of the Hodge-Riemann relations for surfaces.
\end{proof}

Thus, starting from the realization $L\subseteq \CC^E$, one may seek a smooth projective $\CC$-variety equipped with two semi-ample divisor classes $\alpha,\beta$ such that their intersection degrees yield the appropriate combinatorial sequence.
We explain how this is done for the sequences \ref{sseqw} and \ref{sseqw'}.

\subsubsection{Log-concavity of \ref{sseqw}}\label{subsubsec:WL}
We show the slightly stronger statement that the closely related sequence $(\overline w_0, \ldots, \overline w_{r-1})$ defined by
\[
\frac{1}{1+q}\T_\M(1+q,0) = \sum_{i=0}^{r-1} \overline w_i q^{r-1-i}
\]
is log-concave with no internal zeros.
\Cref{thm:poin} states that this sequence is exactly the Betti numbers of the arrangement complement $\PP \mathring L$.
The sought-after projective variety is the {wonderful compactification} of a hyperplane arrangement complement introduced in \cite{dCP95}.

\begin{defn}\label{defn:WL}
The \textbf{wonderful compactification} $W_L$ is the variety obtained by blowing-up $\PP L$ at all the points $\{\PP L_F : \rk_\M(F) = r-1\}$, then by blowing-up all strict transforms of the lines $\{\PP L_F: \rk_\M(F) = r-2\}$, and so forth.  Let $\pi_L: W_L \to \PP L$ be the blow-down map.
\end{defn}

By construction $\pi_L$ is isomorphism on the open loci $\PP \mathring L$.  The boundary $\partial W_L = W_L \setminus \PP \mathring L$ is a simple-normal-crossing divisor on $W_L$ \cite{dCP95}.  For the sought-after divisor classes $\alpha,\beta$ on $W_L$, the wonderful compactification for the Boolean matroid plays a special role.

\medskip
When $\M = \U_{n,n}$, that is, when $L = \CC^E$, the wonderful compactification is known as the \textbf{permutohedral variety} $X_{A_{n-1}}$.  Explicitly, it is obtained from $\PP^{n-1}$ by blowing-up all $n$ coordinate points of $\PP^{n-1}$, then blowing-up all strict transforms of $\binom{n}{2}$ coordinate lines of $\PP^{n-1}$, and so forth.  Let $\pi_1: X_{A_{n-1}} \to \PP^{n-1}$ be the blow-down map. Blowing-down the exceptional divisors in a ``reverse manner'' (see \cite{Huh18b} for a detailed description via toric geometry), one obtains a different blow-down map $\pi_2: X_{A_{n-1}}\to \PP^{n-1}$.  The resulting birational transformation is the Cremona transformation $\operatorname{crem}: \PP^{n-1} \dashrightarrow \PP^{n-1}$ given by $[x_1, \ldots, x_n] \mapsto [\frac{1}{x_1}, \ldots, \frac{1}{x_n}]$ in the projective coordinates.

\medskip
Returning to the case where $\M$ is not necessarily Boolean, we note the following:  Because the arrangement $\mathcal A$ on $\PP L$ is the restriction to $\PP L$ of the coordinate hyperplane arrangement on $\PP^{n-1}$, the universal property of blow-ups implies that $W_L$ is the strict transform of $\PP L \subseteq \PP^{n-1}$ under the blow-up $\pi_1$.  Summarizing, we have a commuting diagram:
\begin{equation}\tag{$\dagger$}\label{fig:WL}
\begin{tikzcd}
&W_L \ar[r, hook]\ar[d, "\pi_L"'] &X_{A_{n-1}} \ar[d, "\pi_1"']\ar[dr, "\pi_2"]\\
&\PP L \ar[r, hook] &\PP^{n-1} \ar[r, dashed, "\operatorname{crem}"]& \PP^{n-1}
\end{tikzcd}
\end{equation}

Let $h$ be the hyperplane class of $\PP^{n-1}$, and define divisor classes $\alpha$ and $\beta$ on $W_L$ to be the restrictions of $\pi_1^*h$ and $\pi_2^*h$, respectively.  Huh--Katz \cite{HK12} showed that
\[
\frac{1}{1+q}T_\M(1+q,0) = \sum_{i=0}^{r-1} \left(\int_{W_L} \alpha^{r-1-i} \beta^{i}\right)q^{r-1-i}.
\]
Since $\alpha$ and $\beta$ are both hyperplane class pullbacks, they are globally-generated, so \Cref{prop:KT} implies the desired log-concavity.

\medskip
How might one think to do this, at least in hindsight?  
Two key steps are as follows.
\begin{enumerate}
\item The commuting diagram \eqref{fig:WL} shows that $\alpha$ is also the pullback to $W_L$ of the hyperplane class in $\PP L$, so we may loosely interpret multiplication by $\alpha$ as restriction to a general hyperplane $H$ in $\PP L$.  As a linear subvariety $H\subset \PP^{n-1}$, this hyperplane $H$ is again a realization of a matroid, known as the \textbf{truncation matroid} $\operatorname{tr}(\M)$ of $\M$.
With well-known properties of characteristic polynomials \cite{Zas87}, it is straightforward to verify that $\frac{1}{1+q}\T_{\operatorname{tr}(\M)}(1+q,0)$ is obtained from $\frac{1}{1+q}\T_{\M}(1+q,0)$ by erasing the constant term and then dividing by $q$.
That is, the sequence $(\overline w_0, \ldots, \overline w_{r-2})$ for $\operatorname{tr}(\M)$ is obtained from that of $\M$ by simply removing the last entry.
\item With the previous step, we now only need compare $\int_{W_L} \beta^{r-1}$ with the constant term $\T_\M(1,0)$.  Let $\PP L^{-1}$ be the closure of the image of $\PP \mathring L$ under the Cremona transformation, often known as the \textbf{reciprocal linear space}.  By the construction of $\beta$, the degree of $\PP L^{-1}$ as a subvariety of $\PP^{n-1}$ equals $\int_{W_L} \beta^{r-1}$.  On the other hand, the degree of $\PP L^{-1}$ also equals $\T_\M(1,0)$.  This last key fact was proven in several contexts \cite{Ter02, PS06, Huh12, HK12}.  A topological approach in \cite{Huh12} is as follows: A result of Dimca and Papadima \cite{DP03} from (complex) Morse theory related the Euler characteristic of a hypersurface complement to the degree of the gradient map.
The hypersurface $\{x_1x_2\cdots x_n = 0\}\subset\PP^{n-1}$ is the coordinate hyperplane arrangement that restricts to the hyperplane arrangement $\mathcal A$ on $\PP L$, and the gradient map of $x_1x_2\cdots x_n$ is exactly the Cremona map.  Combining these facts with \Cref{thm:poin}, one can deduce $\deg \PP L^{-1}= \T_\M(1,0)$.
\end{enumerate}

\begin{exer}
Let $L$ be as in \Cref{eg:U34}.  Verify that $\T_\M(1,0) = 3$, and verify that $\PP L^{-1}$ is a cubic surface known as the Cayley nodal cubic.  This cubic surface has four singular points, with a line through each pair of points; explain where these come from in terms of the wonderful compactification $W_L$.  (Bonus: This cubic surface has three more lines, for the total of nine; where do they come from?)
\end{exer}

\subsubsection{Log-concavity of \ref{sseqw'}}\label{subsubsec:WLN}

The sought-after projective variety is the \textbf{variety of critical points} $\mathfrak X$, formally introduced in \cite{CDFV11} but implicit in previous works related to Varchenko's problem on critical points of master functions on an affine hyperplane arrangement \cite{Var95}.
Here, in order to build upon our previous discussion in \Cref{subsubsec:WL}, we follow \cite[Section 8]{BEST} to describe a smooth birational model of $\mathfrak X$ in terms of the wonderful compactification $W_L$, although it differs slightly from the original description in \cite{Huh15}.

\medskip
Consider the embedding $W_L \hookrightarrow X_{A_{n-1}}$ in the diagram \eqref{fig:WL}.  Let $\mathcal N = \mathcal N_{W_L/X_{A_{n-1}}}$ be the normal bundle, and let $\mathfrak X_L = \PP_{W_L}(\mathcal N^\vee)$ be the projectivization\footnote{Our convention for the projectivization of a vector bundle $\mathcal  E$ on a variety $X$ is that $\PP_X(\mathcal  E) = \operatorname{Proj}_X \operatorname{Sym}^\bullet(\mathcal E^\vee)$, which agrees with \cite{EH16} but is the opposite of \cite{Har77, Laz04b}.} of the conormal bundle with the projection map $p: \mathfrak X_L \to W_L$.  Recall the blow-down map $\pi_L: W_L \to \PP L$.

\medskip
The sought-after divisor classes on $\mathfrak X_L$ are as follows.  Let $\gamma$ be the pullback of the hyperplane class in $\PP L$ via the composition $\mathfrak X_L \overset{p}\to W_L\overset{\pi_L}\to \PP L$.  Let $\delta = c_1(\mathcal O(1))$ be the first Chern class of the line bundle $\mathcal O(1)$ from the construction of $\mathfrak X_L$ as a projectivization of a vector bundle, which turns out to be semi-ample.  One can then translate \cite[Theorem 1.1]{DGS12} to the statement that
\[
\frac{1}{q}\T_\M(q,0) = \sum_{i = 0}^{r-1}\left( \int_{\mathfrak X_L} \gamma^{r-1-i} \delta^{n-r-1+i}\right)q^{r-1-i}.
\]
\Cref{prop:KT} now implies the desired log-concavity.

\medskip
How might one think to do this, at least in hindsight?
For the original formulation of $\mathfrak X$, maximum likelihood problems in algebraic statistics provided a motivation; see \cite{Huh13, HS14}.
For the related construction $\mathfrak X_L$ here, we highlight some key steps.  In either cases, one uses properties of log-tangent bundles and their characteristic classes; see \cite{Alu05} for an introduction to these tools.

\begin{enumerate}
\item By \Cref{thm:poin}, the constant term equals
\[
\left(\frac{1}{q}\T_\M(q,0)\middle)\right|_{q = 0} = \left(\frac{1}{1+q}\T_\M(1+q,0)\middle)\right|_{q = -1} = (-1)^{r-1} \chi_{top}(\PP \mathring L),
\]
the signed Euler characteristic of $\PP \mathring L$.
Euler characteristics satisfy the ``scissors relation'' that $\chi_{top}(X) = \chi_{top}(X\setminus Z) + \chi_{top}(Z)$ for a closed embedding $Z\hookrightarrow X$ of $\CC$-varieties.  More generally, Chern-Schwartz-MacPherson (CSM) classes \cite{Mac74} of $\CC$-varieties are homological objects that respect such scissors relation.  One then uses \Cref{exer:scissors} to show that the degrees of the CSM classes of $\PP \mathring L$ satisfy the same deletion-contraction relation satisfied by the coefficients of $\frac{1}{-q}\T_M(-q,0)$, so that \Cref{thm:delcont} implies that they are the same.

\item Having related the CSM classes to $\frac{1}{-q}\T_M(-q,0)$, we now relate powers of $\delta$ with the CSM classes of $\PP \mathring L$.
The varieties $W_L$ and $X_{A_{n-1}}$ have simple-normal-crossing boundaries $\partial W_L$ and $\partial X_{A_{n-1}} = X_{A_{n-1}} \setminus ((\CC^*)^E/\CC^*)$.
Under the embedding $W_L \hookrightarrow X_{A_{n-1}}$, one can show that $\partial W_L = W_L \cap \partial X_{A_{n-1}}$ scheme-theoretically.
Consequently, one has a short exact sequence (see for instance \cite[Section 9]{EHL})
\[
0 \to \mathcal T_{W_L}(-\log \partial W_L) \to \mathcal T_{X_{A_{n-1}}}(-\log \partial X_{A_{n-1}})|_{W_L} \to \mathcal N \to 0.
\]
Because $X_{A_{n-1}}$ is a toric variety with the dense open torus $(\CC^*)^E/\CC^*$, the log-tangent bundle $T_{X_{A_{n-1}}}(-\log \partial X_{A_{n-1}})$ is trivial \cite[Chapter 8]{CLS11}.  Thus, the normal bundle $\mathcal N$ is globally-generated so that $\delta$ is semi-ample.  Moreover, the Segre classes of $\mathcal N$, which are given by powers of $-\delta$, equal the Chern classes of $\mathcal T_{W_L}(-\log \partial W_L)$, which are the CSM classes of the complement $W_L \setminus \partial W_L = \PP \mathring L$ \cite{Alu99}.
\end{enumerate}

\subsection{The top-heaviness via intersection cohomology}\label{subsec:YL}
We start with a general strategy for establishing top-heaviness, which first appeared in \cite{BE09} to establish top-heaviness for Bruhat intervals.  Suppose we have a (not necessarily smooth) projective $\CC$-variety $X$ with an \textbf{affine stratification}: There is a finite collection $\{U_j\}_{j\in J}$ of locally closed subvarieties of $X$, called the \textbf{strata}, such that $X$ is the disjoint union of $\{U_j\}_{j\in J}$, the closure $\overline{U_j}$ of any strata is again a union of strata, and each strata is isomorphic to $\mathbb C^m$ for some $m$.

\begin{thm}\label{thm:BE} \cite[Theorem 3.1]{BE09}
For $0\leq i \leq d = \dim X$, let $b_i$ be the number of strata of dimension $i$.  Then, we have $b_i \leq b_j$ for all $0\leq i \leq j \leq d-i$.
\end{thm}

\begin{proof}[Sketch of the proof]
Let us assume we have shown that $\dim H^{2i}(X) = b_i$.
If $X$ were smooth, the hard Lefschetz theorem implies that $(b_0, \ldots, b_d)$ is in fact unimodal and symmetric, so we would be done.  Since $X$ may not be smooth, we need to consider the intersection cohomology $IH^\bullet(X)$, for which the hard Lefschetz theorem still holds \cite{GM83, BBD82}.  There is a natural graded map $H^\bullet(X) \to IH^\bullet(X)$, fitting into a commutative diagram
\[
\begin{tikzcd}
&H^{2i}(X) \ar[r]\ar[d,"\cdot c_1(\mathcal L)^{j-i}"'] &IH^{2i}(X) \ar[d, hook, "\cdot c_1(\mathcal L)^{j-i}"]\\
&H^{2j}(X) \ar[r] &IH^{2j}(X)
\end{tikzcd}
\]
for any ample line bundle $\mathcal L$ on $X$, where the injectivity of the right vertical map follows from the validity of the hard Lefschetz property for intersection cohomology.  Thus, if $H^\bullet(X) \to IH^\bullet(X)$ is injective, then the left vertical map is necessarily injective, so we can conclude the desired $b_i \leq b_j$.

The proof of $\dim H^{2i}(X) = b_i$ and the injection $H^\bullet(X) \hookrightarrow IH^\bullet(X)$ both follow from combining a standard long exact sequence of cohomologies and the result of Weber \cite{Web04} (see also \cite[Theorem 2.1]{BE09}) which identified the kernel of $H^k(X) \to IH^k(X)$ in terms of the weight filtration on $H^k(X)$ given by the mixed Hodge structure \cite{Del71}.
\end{proof}

Bj\"orner and Ekedahl \cite{BE09} used \Cref{thm:BE} on Schubert varieties of a generalized flag variety to deduce top-heaviness for Bruhat intervals of a finite crystallographic Coxeter group.
In our case, for a realization $L\subseteq \CC^E$ of a matroid $\M$, we consider its \textbf{matroid Schubert variety} $Y_L$ defined as
\[
Y_L = \text{the closure of $L$ inside $(\PP^1)^E$},
\]
where $\CC^E \subset (\PP^1)^E$ via the identification $\PP^1 = \CC \cup \{\infty\}$.
Note that the identification $\PP^1 = \CC \cup \{\infty\}$ induces an affine stratification of $(\PP^1)^E$ with strata $\{\CC^S \times \{\infty\}^{E\setminus S} : S\subseteq E\}$.
One shows that this stratification restricts to give an affine stratification of $Y_L$ by using the computation of (the Gr\"obner basis for) the defining ideal of $Y_L\subseteq (\PP^1)^E$ given by Ardila and Boocher \cite{AB16}.

\begin{thm}\label{thm:YL} \cite[Theorem 14]{HW17}, \cite[Lemmas 7.5 \& 7.6]{PXY18}
The matroid Schubert variety $Y_L$ admits an affine stratification by the strata $\{U^F : F \text{ a flat of $\M$}\}$ defined by
\[
U^F = Y_L \cap \big(\CC^F \times \{\infty\}^{E\setminus F}\big).
\]
For each flat $F$, the strata $U^F$ is isomorphic to the image of the composition $L\twoheadrightarrow \CC^E \twoheadrightarrow \CC^{F}$, which has dimension $\rk_\M(F)$.
\end{thm}

The top-heaviness of the sequence \ref{sseqW} now follows from Theorems~\ref{thm:YL} and \ref{thm:BE}.

\begin{exer}
Verify \Cref{thm:YL} for a realization of the uniform matroid $\U_{2,3}$.
\end{exer}

For Schubert varieties in generalized flag varieties, their intersection cohomology is closely related to Kazhdan-Lusztig theory.
The terminology ``matroid Schubert variety'' was chosen because of the analogous relation between the intersection cohomology $Y_L$ and Kazhdan-Lusztig theory of matroid developed in \cite{EPW16, PXY18}.

\section{Tropical Hodge theory}\label{sec:tropHodge}

We will begin by explaining how polyhedral fans and matroids give rise to ``cohomology'' rings.  We then discuss two fundamental theorems for these ``cohomology'' rings concerning the validity of the K\"ahler package (\Cref{defn:Kahler}).
We then explain how the log-concavity statements for arbitrary matroids can be deduced from these fundamental theorems.
We will assume some familiarity with the basics of polyhedral geometry.  For a brief introduction see \cite[Chapter 1.2]{Ful93}, and see \cite{Zie95} for a fuller treatment.

\subsection{Chow rings of fans and matroids}

Let $N$ be a lattice, i.e.\ a finitely generated free abelian group $\ZZ^m$.  Let $N^\vee$ denote its dual lattice.  We write $N_\RR = N\otimes \RR$.
Recall that a fan $\Sigma$ in $N_\RR$ is \emph{rational} if each ray $\rho$ in $\Sigma$ equals $\RR_{\geq 0} u$ for some $u\in N$, \emph{simplicial} if every $k$-dimensional cone in $\Sigma$ is generated by $k$ rays, and \emph{pure-dimensional} if every maximal cone has the same dimension.  For each ray $\rho$, let $u_\rho \in N$ be the \emph{primitive ray generator}, i.e.\ the element such that $\rho \cap N = \ZZ_{\geq 0} u_\rho$.  Let $\Sigma(1)$ denote the set of rays of $\Sigma$.  The \emph{support} of $\Sigma$
is denoted $|\Sigma|$.  A fan $\Sigma$ in $N_\RR$ is \emph{complete} if $|\Sigma| = N_\RR$.

\proof[Assumption] All fans we treat will be rational, simplicial, and pure-dimensional, but not necessarily complete.

\begin{defn}
The \textbf{Chow (cohomology) ring} (with real coefficients) of a fan $\Sigma$ in $N_\RR$ is the graded $\RR$-algebra
\[
A^\bullet(\Sigma) = \frac{\RR[x_\rho: \rho\in \Sigma(1)]}{I_\Sigma + J_\Sigma}\]
where $I_\Sigma$ and $J_\Sigma$ are the ideals
\begin{align*}
I_\Sigma &= \Big\langle \prod_{\rho\in S} x_\rho : \textnormal{$S\subseteq \Sigma(1)$ do not form a cone in $\Sigma$}\Big\rangle \quad\text{and}\\
J_\Sigma &= \Big\langle \sum_{\rho\in \Sigma(1)}  m(u_\rho) x_\rho : m \in N^\vee \Big\rangle.
\end{align*}
\end{defn}

It is an exercise to show that the $k$-th graded piece $A^k(\Sigma)$ of $A^\bullet(\Sigma)$ is generated by square-free monomials of degree $k$ in the variables.  In particular, $A^k(\Sigma) = 0$ for all $k$ greater than the dimension $d$ of $\Sigma$, and $A^\bullet(\Sigma) = \bigoplus_{i=0}^d A^i (\Sigma)$ is a finite-dimensional graded real vector space.

\begin{exer}
For two fans $\Sigma$ and $\Sigma'$, show that $A^\bullet(\Sigma\times \Sigma') = A^\bullet(\Sigma) \otimes A^\bullet(\Sigma')$.
\end{exer}

\medskip
Borrowing language from algebraic geometry, let us call a linear combination of the variables $x_\rho$ a \textbf{divisor} and its image in $A^1(\Sigma)$ its \textbf{divisor class} on $\Sigma$.
Because $\Sigma$ is simplicial, a divisor $D = \sum_{\rho \in \Sigma(1)} c_\rho x_\rho$ determines a piecewise-linear function $\varphi_D$ on $|\Sigma|$ by assigning the value $c_\rho$ to each primitive ray generator $u_\rho$.

\begin{defn}
A divisor $D$ on a complete fan $\Sigma$ is \textbf{ample} if the piecewise-linear function $\varphi_D$ is strictly-convex, i.e.\ $\varphi_D(u) + \varphi_D(v) < \varphi_D(u+v)$ for all $u,v\in N_\RR$ not in the same cone of $\Sigma$.  It is \textbf{nef} if only the weak inequality $\leq$ is satisfied.

For a not necessarily complete fan $\Sigma$, a divisor $D$ is \textbf{ample} (resp. \textbf{nef}) if $\varphi_D$ is the restriction of the piecewise-linear function of an ample (resp.\ nef) divisor on a complete fan containing $\Sigma$ as a subfan.
We denote by $\mathcal K(\Sigma) \subset A^1(\Sigma)$ the convex set of the divisor classes of ample divisors on $\Sigma$.\footnote{Technically, our definition of ample/nef divisors on non-complete fans differs from that of \cite{ADH22}, which makes our $\mathcal K(\Sigma)$ a subset of the one in \cite{ADH22}, but will suffice for our discussion.}
We often consider $\mathcal K(\Sigma)$ as a set of graded linear maps $A^\bullet(\Sigma) \to A^{\bullet+1}(\Sigma)$ given by multiplication.
\end{defn}

In terms of toric geometry, the ring $A^\bullet(\Sigma)$ is the Chow cohomology ring of the toric variety $X_\Sigma$ associated to the fan $\Sigma$ \cite{Dan78, Bri96}.  When $\Sigma$ is the normal fan of a simple polytope, the toric variety $X_\Sigma$ is a projective variety with mild (i.e.\ orbifold) singularities, whose ample cone is $\mathcal K(\Sigma)$.  Classical results in algebraic geometry then imply the validity of the K\"ahler package for the ring $A^\bullet(\Sigma)$.
Stanley used this to resolve McMullen's $g$-conjecture on the number of faces of a simple polytope \cite{Sta80}.  Afterwards, McMullen \cite{McM93} gave a purely combinatorial proof of the K\"ahler package that works even for nonrational fans (with no associated toric variety in the background).
In our case, the fans will not in general be complete, so a priori there is no reason to expect any validity of the K\"ahler package, since the associated toric varieties are not in general compact.  The miraculous result is that certain fans from matroids turn out to enjoy the K\"ahler package.

\medskip
For a matroid $\M$ on $E$ of rank $r$, we construct a fan introduced and studied in \cite{Stu02, AK06, Spe08} as tropical geometric analogue of linear spaces.
The fan will be in the real vector space over the lattice $\ZZ^E/\ZZ\be_E$.  For a subset $S\subseteq E$, let us denote by $\overline\be_S$ the image of $\be_S = \sum_{i\in S} \be_i\in \RR^E$ under the quotient map $\RR^E \to \RR^E/\RR\be_E$.

\begin{defn}
The \textbf{Bergman fan} $\Sigma_\M$ of a rank $r$ matroid $\M$ is a pure $(r-1)$-dimensional fan in $\RR^E/\RR\be_E$ consisting of the maximal cones $\RR_{\geq 0}\{\overline\be_{F_1}, \ldots, \overline\be_{F_{r-1}}\}$, one for each maximal chain $\varnothing \subsetneq F_1 \subsetneq F_2 \subsetneq \cdots \subsetneq F_{r-1} \subsetneq E$ of nonempty proper flats of $\M$.
\end{defn}

\begin{eg}
The Bergman fans of $\U_{2,3}$ and $\U_{3,4}$ are depicted below.
\begin{figure}[h]
\begin{tabular}{cc}
\raisebox{20pt}{
\includegraphics[height=30mm]{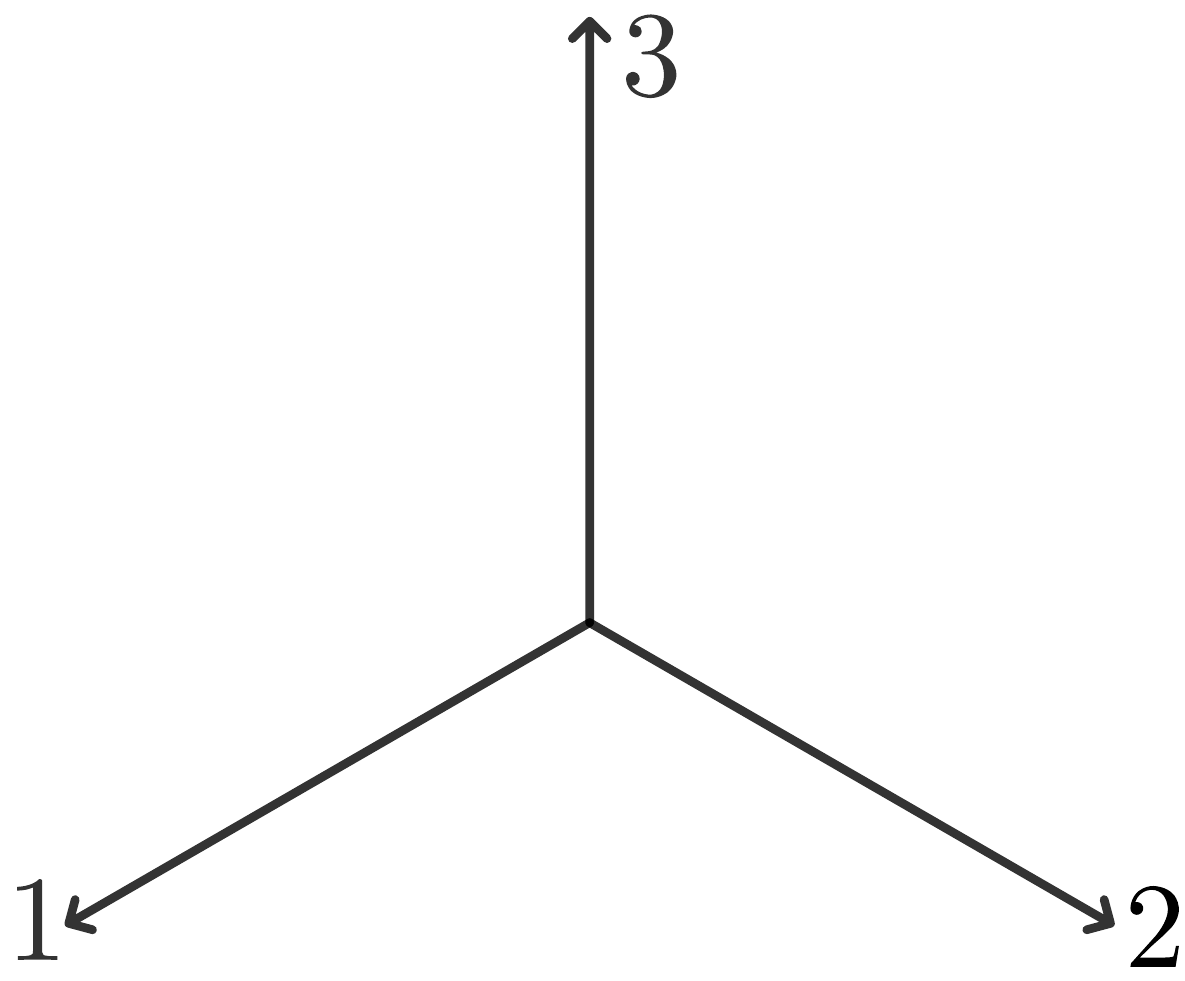}
}
&\qquad\qquad
\includegraphics[height=41.5mm]{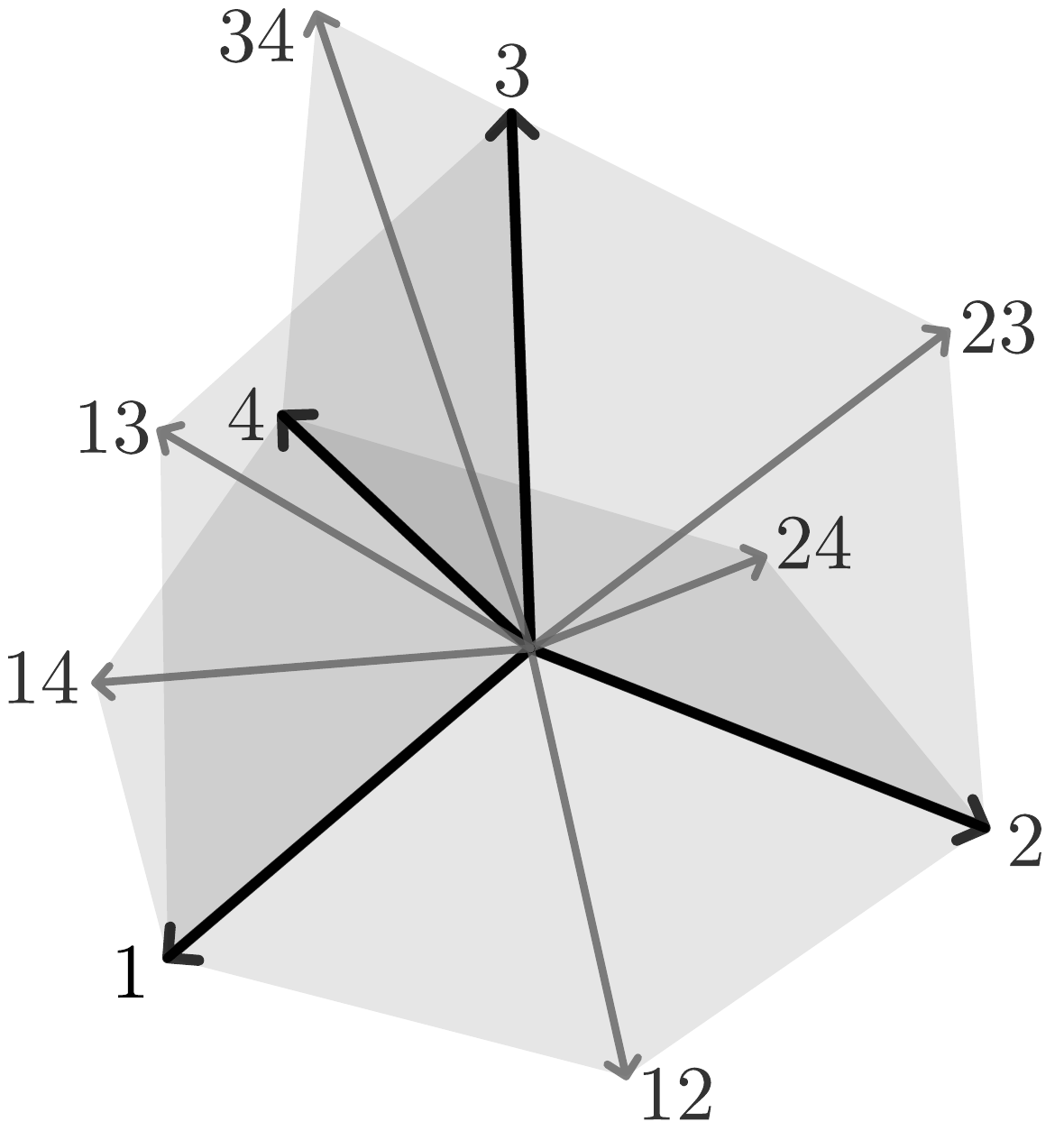}
\end{tabular}
\end{figure}
\end{eg}

\begin{eg}\label{eg:permuto}
For the Boolean matroid $\U_{n,n}$, note that the maximal cones of its Bergman fan $\Sigma_{\U_{n,n}}$ correspond to permutations of the ground set $E$.
The fan $\Sigma_{\U_{n,n}}$ is known as the \textbf{permutohedral fan} or as the \textbf{braid arrangement}, denoted $\Sigma_{A_{n-1}}$.
Note that by construction the Bergman fan of any matroid $\M$ (on ground set $E$) is a subfan of $\Sigma_{A_{n-1}}$.
See \cite{AA} for a survey of remarkable combinatorial properties of the permutohedral fan.
\end{eg}

\begin{defn}
The \textbf{Chow ring $A^\bullet(\M)$ of a matroid} $\M$ is the Chow ring $A^\bullet(\Sigma_\M)$ of its Bergman fan.  Explicitly, it is
\[
A^\bullet(\M) = \bigoplus_{i=0}^{r-1} A^i(\M) = \frac{\RR[x_F: F \text{ a nonempty proper flat of }\M]}{I_\M + J_\M}
\]
where $I_\M$ and $J_\M$ are the ideals
\begin{align*}
I_\M &= \Big\langle x_Fx_{F'} : F\not\subseteq F' \text{ and } F\not\supseteq F'\Big\rangle \quad \text{and}\\
J_\M &= \Big\langle \sum_{F\ni i} x_F - \sum_{G\ni j} x_G: i,j \in E\Big\rangle, \begin{matrix}\text{ where the two sums are over all flats}\\
\text{ containing $i$ and $j$, respectively.}\end{matrix}
\end{align*}
\end{defn}

\begin{rem}\label{rem:chowequiv}
When $\M$ has a realization $L\subseteq \CC^E$, the wonderful compactification $W_L$ defined in \Cref{subsubsec:WL} is the closure of $\PP \mathring L$ inside the toric variety $X_{\Sigma_\M}$.  The resulting pullback map of cohomologies is an isomorphism between the Chow ring $A^\bullet(\M) = A^\bullet(X_{\Sigma_\M})$ of the toric variety of $\Sigma_\M$ and the cohomology ring of $W_L$ \cite{FY04, dCP95}.
This ``Chow equivalence'' is informed by the theory of \emph{tropical compactifications} \cite{Tev07}; see \cite[Chapter 6]{MS15} for an introduction.
While the K\"ahler package for the Chow ring in this realizable case thus follows from classical Hodge theory, \cite[Theorem 5.12]{AHK18} states that the existence of such Chow equivalence is equivalent to the realizability of the matroid.
\end{rem}

We note another fan associated to a matroid whose Chow ring will be used for log-concavity statements.

\begin{defn}
The \textbf{conormal fan} $\Sigma_{\M,\M^\perp}$ of a (loopless and coloopless) matroid $\M$ is a fan in $\RR^E/\RR\be_E \times \RR^E/\RR\be_E$ whose support equals the support of the product $\Sigma_\M \times \Sigma_{\M^\perp}$.
\end{defn}

We omit the precise definition, which involves the intricate and interesting combinatorics of the \emph{bipermutohedron} and \emph{biflags} introduced in \cite{ADH22}.  See \cite[Section 2.8]{ADH22} for its origin story.

\begin{rem}
Just as the Chow ring $A^\bullet(\M)$ of a matroid $\M$ is modeled after the wonderful compactification, the \emph{conormal Chow ring} of a matroid is modeled after the geometry described in \Cref{subsubsec:WLN}.
For instance, \Cref{exer:dual} suggests how the product $\Sigma_\M \times \Sigma_{\M^\perp}$ can serve as a polyhedral model of the projectivized conormal bundle $\mathfrak X_L = \PP_{W_L}(\mathcal N^\vee)$.
\end{rem}

\subsection{Fundamental theorems of tropical Hodge theory}\label{subsec:tropHodgeKahler}

We now discuss two fundamental theorems concerning the K\"ahler package for Chow rings of not necessarily complete fans.  To state them, we need a few more terminologies.

\medskip
For a cone $\sigma$ of a fan $\Sigma$ in $N_\RR$, the \textbf{star} $\operatorname{st}_\sigma(\Sigma)$ is a fan in $N_\RR/\operatorname{span}(\sigma)$ whose cones are the images under the projection $N_\RR \to N_\RR/\operatorname{span}(\sigma)$ of the cones of $\Sigma$ containing $\sigma$.  Geometrically, the toric variety $X_{\operatorname{st}_\sigma(\Sigma)}$ of the star is the closure of the torus-orbit of $X_\Sigma$ corresponding to the cone $\sigma$.

\begin{exer}\label{exer:ann} \
\begin{enumerate}

\item The star of the Bergman fan $\Sigma_\M$ at the ray $\RR_{\geq 0}\be_F$ corresponding to a nonempty proper flat $F$ is isomorphic to the product $\Sigma_{\M|F}\times \Sigma_{\M/F}$.

\item Show that the Chow ring $A^\bullet(\operatorname{st}_\sigma\Sigma)$ of the star is isomorphic to the quotient ring $A^\bullet(\Sigma)/\langle a \in A^\bullet(\Sigma) : a \cdot \prod_{\rho\in\sigma} x_\rho = 0\rangle$.

\end{enumerate}
\end{exer}

\medskip
A \textbf{positive Minkowski weight} on a (pure) $d$-dimensional fan $\Sigma$ is a linear map $\deg: A^d(\Sigma) \to \RR$ such that $\deg(\prod_{\rho\in \sigma}x_\rho) > 0$ for any maximal cone $\sigma$ of $\Sigma$.  
It defines a symmetric bilinear pairing $A^\bullet \times A^{d -\bullet} \to \RR$ by $(x,y) \mapsto \deg(xy)$, which by abuse of notation we also denote $\deg$.

Geometrically, Minkowski weights in general are fundamental objects in tropical intersection theory, serving the role of ``Chow homology classes.''  For their definition and properties, see \cite{FS97, KP08} and \cite[Chapter 6]{MS15}, as well as \cite[Section 5]{AHK18}.  Fans arising in the context of tropical compactifications (\Cref{rem:chowequiv}) provide many examples of positive Minkowski weights.  For instance, the following can be deduced from what is known as the cover-partition property of flats.

\begin{prop}\cite[Proposition 5.2]{AHK18}
For a matroid $\M$ of rank $r$, the assignment $x_{F_1} \cdots x_{F_{r-1}}\mapsto 1$ for any maximal chain $\varnothing \subsetneq F_1 \subsetneq \cdots \subsetneq F_{r-1}\subsetneq E$ of nonempty proper flats of $\M$ gives a well-defined linear map $\deg_\M: A^{r-1}(\M) \to \RR$.
\end{prop}

Similarly, the conormal fan of a matroid also has a natural positive Minkowski weight $\deg_{\M,\M^\perp}$. 
The key notion for the statement of the two fundamental theorems is the notion of Lefschetz-ness of fans, introduced in \cite{ADH22}.

\begin{defn}
A fan $\Sigma$ of dimension $d$ is said to be \textbf{Lefschetz} if the following are satisfied:
\begin{enumerate}
\item $\operatorname{Hom}_\RR(A^d(\Sigma), \RR)$ is spanned by a positive Minkowski weight $\deg$.
\item The triple $(A^\bullet(\Sigma),\deg, \mathcal K(\Sigma))$ satisfies the K\"ahler package (\Cref{defn:Kahler}).
\item For any positive dimensional cone $\sigma$ of $\Sigma$, the star $\operatorname{st}_\sigma(\Sigma)$ is Lefschetz.
\end{enumerate}
\end{defn}

We can now state the two fundamental theorems.

\begin{thm}\label{thm:Bergman} \cite[Theorem 1.4]{AHK18}
The Bergman fan $\Sigma_\M$ of a matroid $\M$ is Lefschetz.
\end{thm}

\begin{thm}\label{thm:tropHodge} \cite[Theorem 1.6]{ADH22}
Let $\Sigma$ and $\Sigma'$ be fans in $N_\RR$ with the same support, and suppose $\mathcal K(\Sigma)$ and $\mathcal K(\Sigma')$ are nonempty.  Then $\Sigma$ is Lefschetz if and only if $\Sigma'$ is Lefschetz.
\end{thm}

The product of two Lefschetz fans is again Lefschetz \cite[Section 7.2]{AHK18}.  Combining this with the two theorems yields the following.

\begin{cor}\label{cor:conormal}
The conormal fan of the matroid is Lefschetz.
\end{cor}

\begin{exer}\label{exer:rank3}
Let $\M$ be a matroid of rank $3$.  Let $\alpha = \sum_{G\ni i} x_G$ for any $i\in E$.  Note that $\alpha\in A^\bullet(\M)$ is independent of the choice of $i$ due to the linear relations $J_\M$.
\begin{enumerate}
\item Show that $\deg_\M(x_F^2) = -1$ for any flat $F$ of rank 2.
\item Show that $\deg_\M\big( \alpha^2 \big) = 1$ for any element $i\in E$.
\item Show that $\deg_\M(\alpha x_F) = 0$ for any element $i\in E$ and a flat $F$ of rank 2.
\item Use these steps to establish the K\"ahler package for $A^\bullet(\M)$ with $\mathcal K = \RR_{\geq 0} \alpha$.
\end{enumerate}
\end{exer}

Let us give a broad overview of the proofs of the two theorems.  
Both employ the following strategy for establishing the K\"ahler package for a graded $\RR$-algebra $A^\bullet$ of ``dimension $d$'' in the sense that $A^\bullet = \bigoplus_{i = 0}^d A^i$.
This general strategy and variations thereof appear in several previous works on the K\"ahler package across varied mathematical fields, such as the works of McMullen \cite{McM93} on simple polytopes, Elias and Williamson \cite{EW14} on Soergel bimodules, and \cite{dCM09} on the topology of algebraic maps.

\begin{enumerate}[label = (\roman*)]
\item It suffices to show the statements of (HL) and (HR) in \Cref{defn:Kahler} in the special case where $L_0 = L_1 = \cdots = L_{d-2i}$, because this ``non-mixed'' version of the K\"ahler package implies the original ``mixed'' version of the K\"ahler package \cite{Cat08} (see also \cite[Theorem 5.20]{ADH22}).
\item Next, one can set up an induction on the ``dimension'' $d$ as follows.   In many situations, the quotient algebra $A^\bullet/\operatorname{ann}(x)$ for certain choices of $x\in A^1$ is again in the family of graded $\RR$-algebras that one is seeking to establish the K\"ahler package for.  Here, $\operatorname{ann}(x)$ denotes the annihilator $\{a \in A^\bullet : ax = 0\}$.  The key observation then is that
\[
\begin{pmatrix}
\text{(HR) in degree $i$ of $A^\bullet/\operatorname{ann}(x)$}\\
\text{for sufficiently many $x\in A^1$}
\end{pmatrix} 
\implies
\begin{pmatrix}
\text{(HL) in degree $i$}\\
\text{for the original $A^\bullet$}
\end{pmatrix} 
\]
(see for instance \cite[Proposition 6.1.6]{BES}).  Because the quotient $A^\bullet/\operatorname{ann}(x)$ has smaller ``dimension'', i.e.\ its $d$-th graded part is zero, one can now proceed by induction on $d$.
\item By the validity of (HL) from the inductive hypothesis, the validity of (HR) for any single element $L\in \mathcal K$ then implies (HR) for all $L\in \mathcal K$.  Thus, the last step is to finish the induction by establishing (HR) for a well-chosen $L\in \mathcal K$.  This is often the most intricate step.
\end{enumerate}

\smallskip
Returning to the case of matroids, step (ii) is provided by \Cref{exer:ann}, which showed that
$A^\bullet(\M)/\operatorname{ann}(x_F) \simeq A^\bullet(\Sigma_{\M|F}\times \Sigma_{\M/F})$
for a nonempty proper flat $F$ of $\M$.
Hence, one can induct on the rank of the matroid.
In the case of \Cref{thm:tropHodge}, step (ii) is essentially built into the definition that the stars of Lefschetz fans are Lefschetz, so that one can induct on the dimension of the fan.

\medskip
In the original proof \cite{AHK18} of \Cref{thm:Bergman}, in order to carry out step (iii), Adiprasito, Huh, and Katz introduced the notion of ``flips,'' inspired by the proof of the K\"ahler package for simple polytopes by McMullen \cite{McM93}.
This process converts the Bergman fan of $\Sigma_\M$ through a sequence of fans until it reaches a fan for which (HR) can be verified easily, and the process is set up such that the validity of (HR) for any one fan in the sequence implies (HR) for all fans in the sequence.
Geometrically, one may interpret the process of flips as a combinatorial abstraction of the process of constructing the wonderful compactification $W_L$ as a sequence of blow-ups (\Cref{defn:WL}).

\medskip
Afterwards, it was recognized that a key property of \emph{semi-small maps} \cite{dCM02} inspires a strategy that can greatly simplify step (iii).
For a map $f: X\to Y$ of smooth projective varieties, the pullbacks to $X$ of ample divisors on $Y$ generally fail (HL) and (HR) in the cohomology ring of $X$, since the pullbacks are generally not ample on $X$ but only nef (i.e.\ is a limit of ample divisors).
However, a characterizing property of a semi-small map is that the pullbacks still satisfy (HL) and (HR).
This inspires the following approach: One can look for a map $\widetilde A^\bullet \to A^\bullet$ of graded algebras, behaving like a pullback along a semi-small map.  If (HR) is known to hold for $\widetilde A^\bullet$, say by induction, step (iii) would follow.

\smallskip
This insight allowed Braden, Huh, Matherne, Proudfoot, and Wang to give a considerably simplified proof of \Cref{thm:Bergman} in \cite{BHMPW22}.  Using that the deletion $\M\setminus e$ of a matroid $\M$ by a non-coloop element $e$ behaves like a semi-small map, they carry out step (iii) by an induction that reduces to the case of Boolean matroids.
This insight on semi-small maps is also essential in the proof of  \Cref{thm:tropHodge} by Ardila, Denham, and Huh.
A key step \cite[Theorem 5.9]{ADH22}, building upon the works \cite{Wlo97, AKMW02}, states that any two fans with the same support can be related by a sequence of edge stellar subdivisions, which are operations on fans that play the role of semi-small maps in toric geometry.

\subsection{Applications of Hodge-Riemann relations in degree 1}

The K\"ahler package gives rise to log-concave sequences in the following way.

\begin{prop}\label{prop:KTcomb}
Suppose $\Sigma$ is a Lefschetz fan of dimension $d$ with a positive Minkowski weight $\deg$.  Then, for any nef divisor classes $\alpha$ and $\beta$, 
\[
\text{the sequence $(a_0, a_1, \ldots, a_d)$ defined by } a_i = \deg(\alpha^{d-i}\beta^i)
\]
is log-concave with no internal zeros.
\end{prop}

\begin{proof}
We may assume that $\alpha,\beta$ are ample, and show that the sequence is strictly positive and log-concave, since a limit of such sequences is necessarily log-concave with no internal zeros.
Strict positivity is then implied by Hodge-Riemann relations (HR) in degree 0.
For log-concavity, (HR) in degree 1, with $L_1\cdots L_{d-2} = \alpha^{d-i-1}\beta^{i-1}$, implies that the symmetric bilinear pairing $A^1(\Sigma) \times A^1(\Sigma) \to \RR$ given by $(x,y) \mapsto \deg(xy \cdot \alpha^{d-i-1}\beta^{i-1})$ has at most one positive eigenvalue.  This implies that the symmetric matrix
\[
\begin{bmatrix}
\deg(\alpha^2 \cdot\alpha^{d-i-1}\beta^{i-1}) & \deg(\alpha\beta\cdot \alpha^{d-i-1}\beta^{i-1}) \\
 \deg(\alpha\beta\cdot\alpha^{d-i-1}\beta^{i-1}) &  \deg(\beta^2\cdot\alpha^{d-i-1}\beta^{i-1})
\end{bmatrix}
\]
cannot be positive definite, but it also cannot be negative definite because all of its entries are positive.  Hence, the determinant of the matrix is non-positive, or equivalently, $a_i^2 \geq a_{i-1}a_{i+1}$.
\end{proof}

Returning to showing log-concavity for a matroid $\M$ of rank $r$, one now searches for appropriate divisor classes on the Bergman fan $\Sigma_\M$ or the conormal fan $\Sigma_{\M,\M^\perp}$.  This step benefits heavily from the geometry of realizable matroids explained in \Cref{sec:real}.

\medskip
The divisor classes for the log-concavity of the sequence \ref{sseqw} come from an involutive symmetry of the permutohedral fan $\Sigma_{A_{n-1}}$ (\Cref{eg:permuto}).
As $-\overline\be_i = \overline\be_{E\setminus i}$, the minus map $x \mapsto -x$ on $\RR^E/\RR\be_E$ gives an involution of $\Sigma_{A_{n-1}}$.  This equips $\Sigma_{A_{n-1}}$ with two distinguished coarsenings to a normal fan of a simplex:
Let $\Sigma_{\Delta}$ (resp.\ be $\Sigma_{\nabla}$) be the fan whose rays are $\{\overline\be_i: i\in E\}$ (resp.\ $\{-\overline\be_i: i\in E\}$) and whose cones are generated by any subsets of the rays with cardinality $\leq n-1$.
The fan $\Sigma_{A_{n-1}}$ is a common refinement of both the fans $\Sigma_\Delta$ and $\Sigma_{\nabla}$.

Let us fix an element $i\in E$, and let $\varphi_\Delta$ be the piecewise-linear function on $\Sigma_\Delta$ given by $\overline\be_i \mapsto 1$ and $\overline\be_j \mapsto 0$ for all $j\neq i$, which is clearly a convex (but not strictly-convex) function.
Restricting this piecewise-linear function to $\Sigma_\M$, it defines a divisor on $\Sigma_\M$ whose divisor class is denoted $\alpha\in A^1(\M)$.  Similarly, starting with the fan $\Sigma_{\nabla}$, we obtain a divisor class $\beta$.  Both $\alpha$ and $\beta$ are nef, and independent of the choice of $i\in E$ we fixed.
Algebraically, we may take $\alpha = \sum_{F\ni i} x_F$ and $\beta = \sum_{F\not\ni i} x_F$ for any choice of $i\in E$.
We have the following.

\begin{prop}\label{prop:alphabeta} \cite[Proposition 5.2]{HK12}, \cite[Proposition 9.5]{AHK18}
With the notations as above, we have
\[
\frac{1}{1+q}\T_\M(1+q,0) = \sum_{i=0}^{r-1} \deg_\M(\alpha^{r-1-i}\beta^i) q^{r-1-i}.
\]
\end{prop}

Combining the proposition with \Cref{thm:Bergman} and \Cref{prop:KTcomb}, we obtain that the coefficients of $\frac{1}{1+q}\T_\M(1+q,0)$ are log-concave with no internal zeros, which implies the same property for the coefficients of $\T_\M(1+q,0)$, i.e.\ the log-concavity of the sequence \ref{sseqw}.
One can also consider other kinds of divisor classes on $\Sigma_\M$ and their values under $\deg_\M$ to study properties of matroids:  See for instance \cite{Eur20, BES, BST, DR22}.

\begin{rem}
In geometric terms, the toric variety of $\Sigma_{A_{n-1}}$ is the permutohedral variety defined in \Cref{subsubsec:WL}.  That the fan $\Sigma_{A_{n-1}}$ coarsens to $\Sigma_\Delta$ and $\Sigma_\nabla$ gives the two blow-down maps $\pi_1, \pi_2: X_{A_{n-1}}\to \PP^{n-1}$, related by the Cremona transformation, which were described in \eqref{fig:WL}.  Thus, the divisor classes $\alpha$ and $\beta$ here and the computations involving them agree with those described in \Cref{subsubsec:WL}.
\end{rem}

\medskip
For the sequence \ref{sseqw'}, recall that the conormal fan of a matroid $\M$ is a fan in $\RR^E/\RR\be_E \times \RR^E/\RR\be_E$ with support equal to the support of $\Sigma_\M \times \Sigma_{\M^\perp}$.  Let
\begin{align*}
&p: \RR^E/\RR\be_E \times \RR^E/\RR\be_E \to \RR^E/\RR\be_E \text{ be the projection to the first factor, and} \\
&s: \RR^E/\RR\be_E \times \RR^E/\RR\be_E \to \RR^E/\RR\be_E \text{ be the addition map $(x,y) \mapsto x+y$}.
\end{align*}
By pulling back the piecewise-linear function $\varphi_\Delta$ on $\RR^E/\RR\be_E$ along these two maps, we obtain two divisor classes $\gamma$ and $\delta$.

\begin{prop}\label{prop:gammadelta}\cite[Theorem 1.2]{ADH22}
With the notations as above, we have
\[
\frac{1}{q}\T_\M(q,0) = \sum_{i = 0}^{r-1} \deg_{\M,\M^\perp}( \gamma^{r-1-i} \delta^{n-r-1+i})q^{r-1-i}.
\]
\end{prop}

Combining the proposition with \Cref{cor:conormal} and \Cref{prop:KTcomb}, we conclude the log-concavity of the sequence \ref{sseqw'}.
The proposition was proved via an intricate combinatorics of biflags in \cite{ADH22} and the Chern-Schwartz-MacPherson classes of matroids \cite{LdMRS20}.
By developing a new framework of \emph{tautological classes of matroids}, Berget, Spink, Tseng and the author proved a formula \cite[Theorem A \& Theorem 9.7]{BEST} that contains both \Cref{prop:alphabeta} and \Cref{prop:gammadelta} as special cases.

\medskip
Note that only a part of the K\"ahler package, the Hodge-Riemann relations in degrees at most 1, was required for concluding log-concavity.  Extracting the essence of the analytic properties behind (HR) in degrees at most 1 leads to the fascinating theory of \textbf{Lorentzian polynomials} \cite{BH20} and (equivalently) \textbf{completely log-concave polynomials} \cite{ALOGV18b}.
One powerful feature of this theory is that it often allows one to reduce to ``dimension 2'' cases, mirroring the feature in classical algebraic geometry that (HR) in degree 1 can be reduced to the case of surfaces.
For instance, by reducing to an analysis of rank 2 matroids, Br\"and\'en and Huh \cite{BH20} and independently Anari, Liu, Oveis-Gharan, and Vinzant \cite{ALOGV18b} proved that the sequence \ref{sseqI} is in fact \emph{ultra-}log-concave, in the sense that
\[
\frac{I_i^2}{\binom{n}{i}^2} \geq \frac{I_{i-1} I_{i+1}}{\binom{n}{i-1}\binom{n}{i+1}} \quad \text{for all $i$,}
\]
which was conjectured by Mason \cite{Mas72}.
Moreover, by analyzing rank 3 matroids (cf.\ \Cref{exer:rank3}), one can use Lorentzian polynomials to give a simplified proof of the log-concavity of the sequence \ref{seqw} \cite{BES, BL}.
We point to \cite[Section 2]{Huh22} for a survey of Lorentzian polynomials and their applications.

\section{Intersection cohomology of a matroid}\label{sec:IH}

We describe the intersection cohomology of a matroid, and its role in showing the top-heaviness of the sequence \ref{sseqW}.
We begin by considering the following graded algebra.

\begin{defn}
For a matroid $\M$ of rank $r$, its \textbf{M\"obius algebra} is a graded $\RR$-algebra $B^\bullet(\M) = \bigoplus_{i=0}^r B^i(\M)$ where $B^i(\M)$ has basis $\{y_F: F \text{ a rank $i$ flat of } \M\}$ for each $0\leq i \leq r$, and multiplication is given by
\[
y_F \cdot y_{F'} = \begin{cases}
y_{F\vee F'} & \text{if $\rk_\M(F) + \rk_\M(F') = \rk_\M(F\vee F')$}\\
0 & \text{otherwise}.
\end{cases}
\]
\end{defn}

A strategy for the top-heaviness is to show that there is an injective linear map $B^i(\M) \to B^{r-i}(\M)$ for every $i\leq r/2$.  The statement of the hard Lefschetz property inspires a candidate for such a map: the multiplication by a power of an element in $B^1(\M)$.
We then immediately face the difficulty that $B^\bullet(\M)$ usually cannot satisfy Poincar\'e duality (PD) or the hard Lefschetz property (HL), since the sequence $(W_0, \ldots, W_r)$ of the dimensions of graded pieces is usually not symmetric.

\smallskip
The intersection cohomology $IH^\bullet(\M)$, introduced in \cite{BHMPWb}, is a graded vector space containing $B^\bullet(\M)$ that ``most efficiently'' amends the failure of (PD) and (HL).
We give a broad outline of their construction and their properties.
The following remark explains some geometric motivation.

\begin{rem}
Recall from \Cref{subsec:YL} the matroid Schubert variety $Y_L$ of a realization $L\subseteq \CC^E$ of a matroid $\M$.  One deduces from \Cref{thm:YL} that the algebra $B^\bullet(\M)$ is the cohomology ring (in even degrees) of the \emph{matroid Schubert variety} $Y_L$ (see \cite[Theorem 14]{HW17}).
The variety $Y_L$ is usually quite singular, which witnesses the failure of (PD) and (HL) for $B^\bullet(\M)$.  Motivated by the proof of \Cref{thm:BE}, one seeks to understand the intersection cohomology $IH^\bullet(Y_L)$, which contains $B^\bullet(\M)$ as a subalgebra.

To do so, let $f: X\to Y_L$ be a resolution of singularities of $Y_L$.  The decomposition theorem of Beilinson, Bernstein, Deligne, and Gabber \cite{BBD82} implies that $H^\bullet(X)$ can be decomposed into a direct sum of $B^\bullet(\M)$-modules, and that $IH^\bullet(Y_L)$ is a direct summand.
In general, computing these decompositions to get a handle on $IH^\bullet(Y_L)$ can be intractible, but for $Y_L$ there is a resolution $f: \widetilde Y_L \to Y_L$ by the \emph{augmented wonderful variety $\widetilde Y_L$of $L$} \cite{BHMPW22} (see also \cite{EHL}), such that its cohomology ring $H^\bullet(\widetilde Y_L)$ and the injection $B^\bullet(\M) \hookrightarrow H^\bullet(\widetilde Y_L)$ have explicit combinatorial descriptions in terms of the matroid $\M$.
\end{rem}

We first find a bigger graded $\RR$-algebra containing $B^\bullet(\M)$ that satisfies the K\"ahler package.
The \textbf{augmented Bergman fan} \cite[Definition 2.4]{BHMPW22} of a matroid $\M$ of rank $r$ is an $r$-dimensional fan in $\RR^E$ closely related to the Bergman fan $\Sigma_\M$.  Its Chow ring has the following explicit description.

\begin{defn}
The \textbf{augmented Chow ring} (with real coefficients) of a matroid $\M$ is the graded $\RR$-algebra
\[
\operatorname{CH}^\bullet(\M) 
=
\frac{\RR[y_i, x_F : i\in E,\ \text{$F$ a (possibly empty) proper flat of $\M$}]}{\widetilde I_\M + \widetilde J_\M}
\]
where $\widetilde I_\M$ and $\widetilde J_\M$ are the ideals
\begin{align*}
\widetilde I_\M &= \Big\langle x_Fx_{F'} : F\not\subseteq F' \text{ and } F\not\supseteq F'\Big\rangle +  \Big\langle y_ix_F: i\notin F\Big\rangle\quad\text{and}\\
\widetilde J_\M &=\Big\langle y_i - \sum_{F\not\ni i} x_F : i\in E\Big\rangle.
\end{align*}
\end{defn}

The augmented Chow ring has the following useful features:
\begin{itemize}
\item The assignment $y_F \mapsto \prod_{i\in F} y_i$ defines an injection $B^\bullet(\M) \hookrightarrow \operatorname{CH}^\bullet(\M)$ of graded $\RR$-algebras \cite[Proposition 2.28]{BHMPW22}.
\item \Cref{thm:Bergman} combined with \Cref{thm:tropHodge} implies that the augmented Bergman fan is Lefschetz, because the support of the augmented Bergman fan of $\M$ can be identified with the support of the usual Bergman fan of the \emph{free co-extension matroid} of $\M$ (see \cite[Section 5.3]{EHL}).   Thus, the Chow ring $\operatorname{CH}^\bullet(\M)$ satisfies the K\"ahler package.
\end{itemize}
Thus, we have found a bigger algebra containing $B^\bullet(\M)$ that satisfies the K\"ahler package.
However, we are not done because $\operatorname{CH}^\bullet(\M)$ is ``too big'': To conclude injectivity properties for $B^\bullet(\M)$, we need the graded linear operators $\mathcal K$ satisfying (HL) on $\operatorname{CH}^\bullet(\M)$ to come from $B^1(\M)$, but this is almost never the case---a positive linear combination of the $y_i$'s usually does not satisfy (HL) on $\operatorname{CH}^\bullet(\M)$.  One instead must consider the following $B^\bullet(\M)$-submodule of $\operatorname{CH}^\bullet(\M)$ that ``most efficiently'' repairs the the failure of (HL) on $B^\bullet(\M)$.

\begin{thm}
Up to isomorphism there is a unique indecomposable $B^\bullet(\M)$-module direct summand of $\operatorname{CH}^\bullet(\M)$ containing $B^\bullet(\M)$.  This direct summand is the \textbf{intersection cohomology} $IH^\bullet(\M)$ of $\M$.
\end{thm}

In fact, the authors of \cite{BHMPWb} establish a canonical decomposition of $\operatorname{CH}^\bullet(\M)$ as a $B^\bullet(\M)$-module, and identify the direct summand $IH^\bullet(\M)$.   This decomposition along with the K\"ahler package for $\operatorname{CH}^\bullet(\M)$ is then fed into a highly intricate version of the general strategy for establishing the K\"ahler package outlined in \Cref{subsec:tropHodgeKahler}, resulting in the following main theorem.

\begin{thm}\cite[Theorem 1.6]{BHMPWb} 
The intersection cohomology $IH^\bullet(\M)$ of a matroid $\M$ satisfies the K\"ahler package with $\mathcal K = \{\sum_{e\in E} c_ey_e : c_e > 0\}$.
\end{thm}

As a corollary, for any positive linear combination $\ell = \sum_{e\in E} c_e y_e$ and $0\leq i \leq j \leq r-i$, we have a commuting diagram
\[
\begin{tikzcd}
&B^{i}(\M) \ar[r,hook]\ar[d,"\cdot\ell^{j-i}"'] &IH^{i}(\M) \ar[d, hook, "\cdot \ell^{j-i}"]\\
&B^{j}(\M) \ar[r,hook] &IH^{j}(\M)
\end{tikzcd}
\]
where the right vertical map is injective by the hard Lefschetz property of $IH^\bullet(\M)$.  The left vertical map is thus injective, and the desired top-heaviness of the sequence \ref{seqW} follows.

\section{Conclusion}
Matroids are combinatorial structures that capture the essence of independence.
There were several conjectures about the behavior of sequences of invariants of a matroid, involving log-concavity or top-heaviness.
June Huh and his collaborators made fundamental contribution to matroid theory \cite{AHK18, ADH22, BHMPWb}, resolving many of these conjectures.  They began by answering the conjectures for realizable matroids using algebraic geometry, a significant step on its own.  Then, with considerable effort, they were able to extract the combinatorial heart, and establish Hodge-theoretic properties for arbitrary, not necessarily realizable, matroids.
This development of the Hodge theory of matroids forms an integral part of the foundation for studying matroids from an algebro-geometric perspective.

\bibliography{Eur_CEB_submit_final.bib}
\bibliographystyle{alpha}

\end{document}